\documentclass{amsart}

\usepackage{graphicx}
\usepackage{diagrams}
\usepackage{color}

\usepackage{amssymb}
\usepackage{amsmath}
\usepackage{amsthm}
\usepackage{amsfonts}
\usepackage{amstext}
\usepackage{amsopn}

\usepackage{pstricks} %needed for psmatrix

\newtheorem{theorem}{Theorem}[section]

\newtheorem{corollary}[theorem]{Corollary}
\newtheorem{lemma}[theorem]{Lemma}

\newtheorem{proposition}[theorem]{Proposition}

\theoremstyle{definition}
\newtheorem{question}[theorem]{Question}
\newtheorem{definition}[theorem]{Definition}
\newtheorem{example}[theorem]{Example}

\newtheorem{remark}[theorem]{Remark}

\theoremstyle{plain}

\newcommand{\R}{\mathbb{R}}
\newcommand{\Z}{\mathbb{Z}}
\newcommand{\Q}{\mathbb{Q}}

\newcommand{\rk}{\operatorname{rank}}

\newcommand{\onto}{\twoheadrightarrow}

\newcommand{\into}{\hookrightarrow}
\newcommand{\sss}{\scriptscriptstyle}

\newcommand{\an}{\ensuremath{A^{\sss {(n)}}_{\sss H}}}
\newcommand{\anp}{\ensuremath{A^{\sss {(n+1)}}_{\sss H}}}
\newcommand{\aone}{\ensuremath{A^{\sss {(1)}}_{\sss H}}}

\newcommand{\bn}{\ensuremath{B^{\sss {(n)}}_{\sss H}}}
\newcommand{\bnp}{\ensuremath{B^{\sss {(n+1)}}_{\sss H}}}

\newcommand{\gn}{\ensuremath{G^{\sss {(n)}}_{\sss H}}}
\newcommand{\gnp}{\ensuremath{G^{\sss {(n+1)}}_{\sss H}}}
\newcommand{\gone}{\ensuremath{G^{\sss {(1)}}_{\sss H}}}

\newcommand{\ak}{\ensuremath{A^{\sss {(k)}}_{\sss H}}}
\newcommand{\akp}{\ensuremath{A^{\sss {(k+1)}}_{\sss H}}}

\newcommand{\bk}{\ensuremath{B^{\sss {(k)}}_{\sss H}}}
\newcommand{\bkp}{\ensuremath{B^{\sss {(k+1)}}_{\sss H}}}

\newcommand{\rhot}{\ensuremath{\rho}}
\newcommand{\st}{\ensuremath{\sigma^{\sss {(2)}}}}
\newcommand{\rhon}{\rho_n}
\newcommand{\wh}{\widehat}

\def\G{\Gamma}
\def\th{\theta}

\def\SB{\mathcal B}
\def\lra{\longrightarrow}

\def\bbr{{\mathbb R}}
\def\bbz{{\mathbb Z}}
\newcommand{\aut}{\operatorname{Aut}}

\theoremstyle{definition}

\theoremstyle{plain}

\usepackage{amscd}

\begin{document}

\title[Homology Cobordism Invariants]{Homology Cobordism
Invariants and the Cochran-Orr-Teichner Filtration of the Link Concordance Group}

\author[Shelly L. Harvey]{Shelly L. Harvey$^\dag$}
%\author{Shelly L. Harvey}
\address{Department of Mathematics, Rice University\\
Houston, TX, 77005}
\email{shelly@math.rice.edu}

\thanks{$^{\dag}$The author was partially supported by an NSF
Postdoctoral Fellowship, NSF DMS-0539044 and a Sloan Research Fellowship}
\date{}

\begin{abstract} For any group G, we define a new characteristic series related to the derived series, that we call the 
torsion-free derived series of G. Using this series and the Cheeger-Gromov $\rho$-invariant, we obtain new 
real-valued homology cobordism invariants $\rho_{n}$ for closed $(4k-1)$-dimensional manifolds.
For $3$-dimensional manifolds, we show that $\{\rho_{n}| n \in \mathbb{N}\}$ is a linearly independent set and for each $n \geq 0$, the image of $\rho_{n}$ is an 
infinitely generated and dense subset of $\R$. 

In their seminal work on knot concordance, T. Cochran, K. Orr, and P. Teichner define a filtration 
$\mathcal{F}_{(n)}^{m}$ 
of the $m$-component (string) link concordance group, called the $(n)$-solvable filtration. They also define a
grope filtration $\mathcal{G}_{n}^{m}$. 
We show that $\rho_{n}$ vanishes for $(n+1)$-solvable links.  Using this, and the non-triviality
of $\rho_{n}$, we show that for each $m \geq 2$, the successive quotients of
the $(n)$-solvable filtration of the link concordance group contain an infinitely generated
subgroup. 
We also establish a similar result for the 
grope filtration.  We remark that for knots ($m=1$), the successive quotients of the $(n)$-solvable filtration are known to be infinite.  
However, for knots, it is unknown if these quotients have infinite rank when $n\geq 3$.
\end{abstract}

\maketitle

\section{Introduction}\label{section:intro}

The main objective of this paper is to investigate the set of $3$-dimensional manifolds up to \emph{homology cobordism}.  To do this, we define, 
for each $n \in \mathbb{N}$, an invariant of  $(4k-1)$-dimensional manifolds ($k\geq 1$) that we
call $\rho_{n}$.  Loosely speaking, $\rho_{n}(M)$ is defined as a 
signature defect closely associated to a term of 
the \emph{torsion-free derived series} of $\pi_{1}(M)$ and (for smooth manifolds) can be interpreted as the Cheeger-Gromov invariant of $M$ associated to the $n^{th}$ torsion-free derived regular cover of $M$.  
Using a derived version of Stallings' Theorem \cite{CH}, we show that $\rho_{n}$ is an invariant of homology cobordism whereas the
Cheeger-Gromov $\rho$ invariant associated to an arbitrary cover is only a priori a homeomorphism
invariant.
\newtheorem*{homthm}{Theorem~\ref{thm:homcob}}
\begin{homthm}If $M^{4k-1}_1$ is rationally homology cobordant to
$M^{4k-1}_2$ $(k \geq1)$ then $\rhon(M_1)=\rhon(M_2)$.
\end{homthm}  

To define the invariant $\rho_{n}$, we first define a new characteristic series $\{G_{H}^{(n)}\}$ of a group $G$ (Section~\ref{section:tfds}) closely related to the derived series, that 
we call the torsion-free derived series and establish its basic properties.   One should view the torsion-free derived series $\{G_{H}^{(n)}\}$ of a
group $G$ as a series that is closely related to the derived series but whose successive quotients $G_{H}^{(n)}/G_{H}^{(n+1)}$ are torsion-free as $\Z[G/G_{H}^{(n)}]$-modules.
This can be compared to the rational derived series $G_{r}^{(n)}$ (studied in \cite{Ha1}) which is a
series whose successive quotients are torsion-free as abelian groups.

Let $M$ be a $(4k-1)$-dimensional manifold and let $G=\pi_{1}(M)$.  $\rho_{n}(M)$ is defined as follows (see Section~\ref{section:def_rhon} for more details).
Recall that for any group $\Lambda$, the $L^{2}$-signature $\st_{\Lambda}$ is a  real-valued homomorphism on the Witt group of hermitian forms on finitely generated projective 
$\mathcal{U}\Lambda$-modules, where $\mathcal{U}\Lambda$ is the algebra of unbounded operators associated to the von Neumann algebra of $\Lambda$.  Following Hausmann~\cite{Haus}, we show
 (see Lemma~\ref{lem:sbord} and 
Corollary~\ref{cor:sbord}) that for every coefficient system $\phi:G \rightarrow G/G^{\sss (n+1)}_{\sss H}$, there exists a positive integer $r$ and a $4k$-dimensional manifold $W$ such 
that the pair $(rM,r\phi)$ is \emph{stably nullbordant} via $(W,\Phi:\pi_{1}(W)\rightarrow \Lambda)$.  Let $h_{W}$ be the intersection form associated to the regular $\Lambda$-covering space of $W$ and
 let $\sigma$ be the ordinary signature function.  
We show in Lemma~\ref{lem:independence} that $\frac{1}{r}( \st_{\Lambda}(h_{W})-\sigma(W))$ is independent of the stable nullbordism $(W,\Phi)$ 
and define $$\rho_n (M)=\frac{1}{r}( \st_{\Lambda}(h_{W})-\sigma(W)).$$

More generally, we define $\rho_{\Gamma}$ for any coefficient system $\phi:\pi_{1}(M) \rightarrow \Gamma$  (note that one can
also define $\rho_{\Gamma}$ via the Cheeger-Gromov construction for smooth manifolds). Hence one 
could also study $\rho_{(n)}(M)$ (respectively $\rho^{lcs}_{n}(M)$), the canonically defined $\rho$-invariant associated to $\Gamma=\pi_{1}(M)/\pi_{1}(M)^{(n+1)}$ (respectively $\Gamma=\pi_{1}(M)/\pi_{1}(M)^{lcs}_{n+1}$) 
where $G^{(n)}$ (respectively $G^{lcs}_{n}$) is the $n^{th}$ term of the derived series (respectively lower central series) of $G$.  Despite the fact that the torsion-free derived series may seem a bit unwieldy, 
we focus on $\rho_{n}$ (rather than $\rho_{(n)}$ or $\rho_{n}^{lcs}$) since it gives an invariant of rational homology cobordism and provides
new information about the structure of the \emph{$(n)$-solvable} and \emph{grope filtrations}
 of the (string) link concordance group (see below and Section~\ref{section:filtrations}).  
 By contrast, $\rho_{(n)}$ is only a homeomorphism invariant.  
One can use Stallings' Theorem \cite{St} and follow through the proof of Theorem  \ref{thm:homcob} to show that $\rho_{n}^{lcs}$ is a homology cobordism invariant.  However, we choose
to use $\rho_{n}$ in our work since it is more directly related to the $(n)$-solvable and grope filtrations of the link concordance group than $\rho_{n}^{lcs}$.

We show that the $\rho_{n}$ are highly non-trivial and independent for $3$-manifolds (Section~\ref{section:non-triviality}).  To accomplish this, we construct an infinite family of examples of $3$-manifolds
 $\{M(\eta,K)\}$ that are constructed by a method known as \emph{genetic infection}. More specifically, to construct $M(\eta,K)$ we start with a  $3$-manifold $M$ and \emph{infect} $M$ by a knot $K$ along a curve
  $\eta$ in the $n^{th}$ term of the derived series of the fundamental group of $M$.  
We prove that $\rho_{i}(M(\eta,K))$ depends only on $n$ and $\rho_{0}(K)$, where $\rho_{0}(K)$ is the integral of the Levine-Tristam signatures of $K$.

\newtheorem*{calc}{Theorem~\ref{diff_rho}}
\begin{calc}
Let $M$ be a compact, orientable manifold, $\eta$ an embedded curve in $M$, $K$ a knot in $S^{3}$, and $P=\pi_1(M)$.  If $\eta \in P_H^{(n)} - P_H^{(n+1)}$ for some
$n\geq 0$ then 
$$\rho_i(M(\eta,K))-\rho_i(M)= 
\left\{
\begin{array}{ll}
    0 & 0 \leq i\leq  n-1;\\
    \rho_0(K) & i \geq n.\\
\end{array}
\right.    
$$   
\end{calc}

Let $\mathcal{H}^3_\Q$ be the set of 
$\Q$-homology cobordism classes of closed, oriented $3$-dimensional manifolds.  Using the set of examples $\{M(\eta,K)\}$ with varying $K$ we establish the following theorem.

\newtheorem*{thmgen}{Theorem~\ref{gen_dense}}
\begin{thmgen}
The image
of $\rho_n : \mathcal{H}^3_\Q \rightarrow \R$ is (1) dense in $\R$
and (2) an infinitely generated subgroup of $\R$ .
\end{thmgen}

\noindent Moreover, in Theorem~\ref{thm:indep}, we show that $\{\rho_{n}\}$ is a linearly independent subset of the vector space  of functions on $\mathcal{H}^3_\Q$.

We remark that S. Chang and S. Weinberger \cite{CW} use a similar type of signature defect, one associated to the universal cover of a manifold, to define a \emph{homeomorphism} invariant $\tau_{(2)}$ of a $(4k-1)$-dimensional manifold ($k \geq 1$).
Using $\tau_{(2)}$ they show that if $M$ is a $(4k-1)$-dimensional (smooth) manifold with $k\geq 2$ and $\pi_{1}(M)$ is not torsion-free
then there are infinitely many (smooth) manifolds homotopy equivalent to $M$ but not homeomorphic to $M$. 

\vspace{5pt}
For the rest of the paper, we turn our attention to the study of link concordance.  Recall 
 that if two links $L_{1}$ and $L_{2}$ in $S^{3}$ are concordant then $M_{L_{1}}$ and $M_{L_{2}}$ are homology cobordant where $M_{L}$ is the zero surgery on $L$.  
 We define $\rho_{n}(L)=\rho_{n}(M_{L})$
for a link $L$ in $S^{3}$.  Hence, by Theorem~\ref{thm:homcob}, $\rho_{n}$ is a link concordance invariant. 

In \cite{COT}, T. Cochran, K. Orr and P. Teichner defined the $(n)$-solvable (and grope) filtration of the knot concordance group $\mathcal{C}$.  
 In \cite{COT} and their two subsequent papers, \cite{COT1,CT}, they showed that the 
quotients $\mathcal{F}_{(n)}/\mathcal{F}_{(n+1)}$ of the $(n)$-solvable filtration of the knot concordance group are non-trivial.  In particular, they showed that 
for $n=1,2$,  $\mathcal{F}_{(n)}/\mathcal{F}_{(n+1)}$ has infinite rank and for all $n\geq 3$, the quotient has rank at least $1$.  It is still unknown if any of the quotients 
is infinitely generated for $n\geq 3$.   In the current paper, we investigate the $(n)$-solvable (and grope) filtration $\mathcal{F}_{(n)}^{m}$
of the string link concordance group $\mathcal{C}(m)$ 
and the subgroup generated by boundary links $\mathcal{B}(m)$ for links
with $m\geq 2$ components.  Since connected sum is not a well-defined operation for links, it is necessary to use string links to obtain a group structure.  Using $\rho_{n}$
we show that, for $m\geq 2$, each of the successive quotients of the $(n)$-solvable filtration of the boundary string link concordance group 
$\mathcal{BF}_{(n)}^{m}$ is infinitely generated.  

\newtheorem*{mainthm}{Theorem~\ref{main1}}
\begin{mainthm}For each $n\geq 0$ and $m \geq 2$, the abelianization of $\SB\mathcal{F}_{(n)}^m/\SB\mathcal{F}_{(n+1)}^m$
has infinite rank.  In particular, for $m\geq 2$, $\SB\mathcal{F}_{(n)}^m/\SB\mathcal{F}_{(n+1)}^m$ is an infinitely generated subgroup of $\mathcal{F}_{(n)}^m/\mathcal{F}_{(n+1)}^m$.
\end{mainthm}

\noindent We note the previous theorem holds ``modulo local knotting'' (see Corollary~\ref{cor:localknotting}), hence this result cannot be obtained using the work of Cochran-Orr-Teichner on knots.
We also prove a similar statement for the grope filtrations $\mathcal{BG}_{n}^{m}$ and $\mathcal{G}_{n}^{m}$ of the boundary and string link concordance groups respectively.

\newtheorem*{grpthm}{Theorem~\ref{main2}}
 \begin{grpthm}For each $n\geq 1$ and $m \geq 2$, the abelianization of 
$\SB\mathcal{G}_{n}^m/\SB\mathcal{G}_{n+2}^m$
has infinite rank.  Hence $\SB\mathcal{G}_{n}^m/\SB\mathcal{G}_{n+2}^m$ is an infinitely generated subgroup of $\mathcal{G}_{n}^m/\mathcal{G}_{n+2}^m$.
\end{grpthm} 

\noindent We also prove that the abelianization of $\SB\mathcal{G}_{n}^m/\SB\mathcal{G}_{n+1}^m$ has non-zero rank 
for $n\geq 2$ and $m\geq 2$ in Proposition~\ref{prop:614}. We conjecture these quotients groups are in fact infinitely genererated.

To prove Theorem~\ref{main1}, we first show that $\rho_{n}$ 
is additive when restricted to $\mathcal{B}(m)$, the subgroup of $\mathcal{C}(m)$ consisting of  $m$ component boundary string links.
We note that $\rho_{n}$ is not additive on $\mathcal{C}(m)$ itself.

\newtheorem*{prop1}{Proposition~\ref{cor:hom}}
\begin{prop1}For each $n\geq 0$ and $m\geq 1$, $\rho_{n}: \mathcal{B}(m) \rightarrow \R$ is a homomorphism.
\end{prop1}

Next, we show that $(n+1)$-solvable links have vanishing $\rho_{n}$.  
Thus, for each $n\geq 0$, $\rho_{n}$ is a homomorphism from 
$\SB\mathcal{F}_{(n)}^m/\SB\mathcal{F}_{(n+1)}^m$ to $\mathbb{R}$.
\newtheorem*{thm:rho0}{Theorem~\ref{rho0}}
\begin{thm:rho0}If a $3$-manifold $M$ is $(n)$-solvable then
 for each $(n)$-solution $W$ and $k \leq n$, the inclusion $i:M \rightarrow W$ induces monomorphisms
 \begin{equation}i_{\ast}:  H_{1}( M; \mathcal{K}( \pi_{1}(W)/\pi_{1}(W)^{\sss (k)}_{\sss H}) )\hookrightarrow H_{1}( W; \mathcal{K}( \pi_{1}(W)/\pi_{1}(W)^{\sss (k)}_{\sss H}) )\end{equation}
  and
\begin{equation}i_{\ast}: \frac{\pi_{1}(M)}{\pi_{1}(M)^{(k+1)}_{H}} \hookrightarrow \frac{\pi_{1}(W)}{\pi_{1}(W)^{(k+1)}_{H}};\end{equation}
and  \begin{equation}\rho_{k-1}(M)=0.\end{equation}  Thus, if $L \in \mathcal{F}_{(n)}$ then $\rho_{k-1}(L)=0$ for $k\leq n$.
\end{thm:rho0}

To complete the proof of Theorem~\ref{main1}, we construct a collection of links that are $(n)$-solvable and have independent $\rho_{n}$.  To do this we perform genetic infection on the $m$-component trivial link using some knot $K$  
along some carefully chosen curve $\eta$ in $F^{(n)}$ where $F$ is the fundamental group of the trivial link.  This produces a collection of boundary links $\{L(\eta,K)\}$ that are $(n)$-solvable and 
such that the image of $\rho_{n}$ restricted to $\{L(\eta,K)\}$ is infinitely generated.  

We remark that the invariant $\rho_{n}$ is related to certain ``finite'' concordance invariants of boundary links associated to $p$-groups.  Suppose $L$ is a boundary link with $m$ components.  Then there is a surjective map $\pi:G \rightarrow F$ where  $G=\pi_{1}(M_{L})$ and $F$ is the free group on $m$ generators.
  By Theorem~4.1 of \cite{CH} (respectively Stallings' Theorem \cite{St}), 
$G/G^{\sss (n)}_{\sss H}\cong F/F^{(n)}$ (respectively $G/G_{n} \cong F/F_{n}$).
Moreover, $M_{L}$ is the boundary of a $4$-dimensional manifold $W$ over $F$.  
Since $F/F^{(n)}$ and $F/F_{n}$ are residually finite p-groups, by work of W. L\"{u}ck and T. Schick, both $\rho_{n}$ and $\rho_{l}^{lcs}$ can be approximated by 
signatures of finite covers of $W$ where the covering groups are $p$-groups.  These finite p-group signatures are closely related to the concordance invariants of boundary links studied by 
S. Friedl \cite{F} and J. C. Cha and K. H. Ko \cite{CK}.

\vspace{12pt}
We finish this paper by mentioning some applications to boundary link concordance in Section~\ref{sec:last}.  In particular, we show $\rho_{k}$ gives a homomorphism from certain gamma groups (modulo automorphisms of the free group) to $\R$, generalizing work of Cappell and Shaneson.
Here, $\SB(n,m)$ is the group of concordance classes of $m$ component, $n$-dimensional boundary disk links in $D^{n+2}$.

\newtheorem*{prop2}{Proposition~\ref{csprop}}
\begin{prop2} For each $n \equiv 1 \mod 4$ with $n>1$, and each $k\geq 0$, there is an induced homomorphism
$$
\tilde{\rho}_k: \widetilde\G_{n+3}(\bbz F\to\bbz)/\aut F\lra\bbr
$$
that factors through $\SB(n,m)$.
\end{prop2}

\section{The Torsion-Free Derived Series}\label{section:tfds}

To define $\rhon$, we must first 
introduce a new characteristic series of a group called the 
torsion-free derived series.  In this section, we will define
the torsion-free derived series and establish some its basic properties.

If $G$ is a group then $G/G^{\sss {(1)}}$ is an abelian group but may
have
$\Z$-torsion. If one would like to avoid $\Z$-torsion then, in
direct
analogy to the rational lower-central series, one can define
$G^{\sss (1)}_r=\{x \in G \mid x^k \in [G,G] \text{ for some }k\neq0
\}$, which is slightly larger than $G^{\sss (1)}$,
so that $G/G^{\sss {(1)}}_r$ is $\Z$-torsion-free. Proceeding in this
way,
defining $G^{\sss (n)}_r$ to be the radical of $[G^{\sss
(n-1)}_r,G^{\sss (n-1)}_r]$, leads to
what has been called the \textbf{rational derived series} of $G$
\cite{Ha1} \cite{C} \cite{CT}. This is the most rapidly
descending series for which the quotients of successive terms are
$\Z$-torsion-free abelian groups. Note that if $N$ is a normal
subgroup of $G$ then
$N/[N,N]$ is not only an abelian group but it is also a right module
over $\Z[G/N]$, where the action is induced from
conjugation in $G$ ($[x]g=[g^{-1}xg]$). To define the
torsion-free derived series, we seek to eliminate
torsion ``in the module sense'' from the successive quotients.  We
define the \textbf{torsion-free derived series} $\gn$ of $G$
 as follows.  First, set $G^{\sss (0)}_{\sss H}=G$. For $n\geq 0$,
suppose inductively that
$\gn$ has been defined and is normal in $G$ (we will show that $\gn$ is normal in $G$ below). Let $T_n$ be the subgroup
of ${\gn}/{[\gn,\gn]}$ consisting of $\Z[G/\gn]$-torsion elements,
i.e. the elements $[x]$ for which there exists some non-zero
$\gamma \in \Z[G/\gn]$, such that $[x] \gamma=0$. (In fact, since it
will be
(inductively) shown below that $\Z [G/\gn]$ is an Ore Domain,
$T_n$ is a submodule). Now consider the
epimorphism of groups:
$$
\gn \xrightarrow{\pi_n} \frac{\gn}{[\gn,\gn]}
$$
and define $\gnp$ to be the inverse image of $T_n$ under
$\pi_n$. Then $\gnp$ is, by definition, a normal subgroup
of $\gn$ that contains $[\gn,\gn]$. It follows inductively that
$\gn$ contains $\gnp$ (and $G^{\sss {(n+1)}}_r$). Moreover,
since $\left.{\gn}\right/{\gnp}$ is the quotient of the module
$\left.{\gn}\right/{[\gn,\gn]}$
by its 
torsion submodule, it is a $\Z[G/\gn]$ torsion-free module
\cite[Lemma 3.4]{Ste}.
Hence the successive quotients of the
torsion-free derived subgroups are torsion-free modules over
the appropriate rings. We define
$G^{_{(\omega)}}_H=\bigcap_{n<\omega}\gn$ as usual. 

We now establish some elementary properties of the torsion-free derived
series of a group.  Recall that a group is \textbf{poly-(torsion-free
abelian)} (often abbreviated \textbf{PTFA})
if it has a finite subnormal series whose successive quotients are
torsion-free abelian groups.  Such a group is solvable, torsion free,
and 
locally indicable \cite[Proposition 1.9]{Str}.  If $G$ is PTFA
then $\Z G$ is an Ore domain and hence admits a classical (right) 
ring of quotients $\mathcal{K}G$, into which $\Z G$ embeds
\cite[pp. 591--592]{P}.   Hence any finitely generated
(right) module $M$ over $\Z[G/\gn]$ has a well-defined rank
that is defined to be the rank of the vector space $M
\otimes_{\Z[G/\gn]}\mathcal{K}(G/\gn)$
\cite[p. 48]{Cohn}. Alternatively the rank can be defined to be the
maximal integer $m$ such that $M$ contains a submodule isomorphic
to $(\Z[G/\gn])^m$.      

\begin{proposition}
For each $0 \leq n < \omega$,
$\gn$ is a normal subgroup of $G$ and $G/\gn$ is a poly-(torsion-free
abelian) group.  Consequently, $\Z[G/\gn]$
is an Ore domain.
\end{proposition}
\begin{proof} We prove this by induction on $n$.  The statement is clear for $n=0$.  Assume $G_{\sss H}^{(n)}$ is a normal subgroup of $G$ and $G/G_{\sss H}^{(n)}$ is PTFA.
Let $x\in G^{(n+1)}_{\sss H}$ and $g\in G$.  Then $x$ and $g^{-1}xg$ lie in $G^{(n)}_{\sss H}$ by assumption.  By definition of the right module structure on 
$G_{\sss H}^{(n)}/[G_{\sss H}^{(n)},G_{\sss H}^{(n)}]$, $\pi_{n}(g^{-1}xg)=\pi_{n}(x)g$.
Since $x\in G^{(n+1)}_{\sss H}$, $\pi_{n}(x)$ is torsion.  To show that $g^{-1}xg\in G^{(n+1)}_{H}$, it suffices to show that $\pi_{n}(x)g$ is torsion.  
Recall that the set of torsion elements of any module over an Ore domain is known to be a submodule \cite[p. 57]{Ste}. Since 
$\Z[G/G_{\sss H}^{(n)}]$ is an Ore domain, it follows that the set of torsion elements in $G_{\sss H}^{(n)}/[G_{\sss H}^{(n)},G_{\sss H}^{(n)}]$ is a submodule. Thus, $\pi_{n}(x)g$ is torsion and hence $G^{\sss (n+1)}_{\sss H}$ is normal in $G$.

Consider the normal series for $G/G^{\sss (n+1)}_{\sss H}$:
\[1=\frac{G^{\sss (n+1)}_{\sss H}}{G^{\sss (n+1)}_{\sss H}} \triangleleft \frac{G^{\sss (n)}_{\sss H}}{G^{\sss (n+1)}_{\sss H}} \triangleleft \cdots \triangleleft \frac{G^{\sss (1)}_{\sss H}}{G^{\sss (n+1)}_{\sss H}} \triangleleft \frac{G}{G^{\sss (n+1)}_{\sss H}}.\]
Since the successive quotients of the above series are torsion-free abelian groups, $G/G^{\sss (n+1)}_{\sss H}$ is PTFA.
\end{proof}

For convenience, we will often write $G/\gn$ as $G_n$ for any group
$G$ (not to be confused with the terms of the 
lower central series of $G$ which we will denote by $G_{n}^{\text{lcs}}$ in this paper). 

We remark that the torsion-free derived subgroups are characteristic
subgroups but they are not \emph{totally invariant}.
That is, an arbitrary homomorphism $\phi: A \rightarrow B$ need not
send $\an$ to $\bn$.  To see this, let
$A = <x,y,z | [z,[x,y]]>$, $B=<x,y>$ and $\phi : A \rightarrow B$ be
defined by $\phi(x)=x$, $\phi(y)=y$ and $\phi(z)=1$.  
Then $[x,y] \in A^{\sss {(2)}}_{\sss H}$ since $[x,y]$ is
$(z_\ast-1)$-torsion in $\aone/[\aone,\aone]$  where $z_\ast = [z]
\in A/A^{_{(1)}}_{^H}$
 but $\phi([x,y])=[x,y] \not\in B^{\sss (2)}=B^{\sss {(2)}}_{\sss H}$
(see Proposition~\ref{torsionfree}).
\vspace{10pt}

\begin{proposition}\label{functoriality} If $\phi :A \to B$ induces a
monomorphism $\phi :A/\an \into B/\bn$, then $\phi
(\anp)\subset \bnp$ and hence $\phi$ induces a homomorphism
$\phi :A/\anp \to B/\bnp$.
\end{proposition}

\begin{proof}Note that the
hypothesis implies that $\phi$ induces a ring monomorphism
$\tilde{\phi}:\Z[A/\an]\to\Z[B/\bn]$. Suppose
that $x\in \anp$. Consider the diagram below.

\begin{diagram}
\an & \rTo^{\pi_{\sss A}} &\frac{\an}{\anp} \\
\dTo_{\phi} & & \dTo_{\bar{\phi}} \\
\bn & \rTo^{\pi_{\sss B}} & \frac{\bn}{\bnp}
\end{diagram}

\noindent By definition, $\pi_{\sss A}(x)$ is torsion. That is, there
is some non-zero
$\gamma\in\Z[A/\an]$ such that $\pi_{\sss A}(x)\gamma =0$. It is
easy to
check that $\bar{\phi}$ is a homomorphism of right
$\Z[A/\an]$-modules using the module structure induced on
${\bn}/{\bnp}$ by $\tilde{\phi}$ (since $\phi (a^{-1}xa)=\phi(a)^{-1}
\phi(x) \phi(a)$). Thus $\bar{\phi}(\pi_{\sss A}(x))\tilde\phi
(\gamma)=0$. Since $\tilde{\phi}$ is injective,
$\bar{\phi}(\pi_{\sss A}(x))$ is a $\Z[B/\bn]$-torsion element. But
$\bar{\phi}(\pi_{\sss A}(x))=\pi_{\sss B}(\phi (x))$, showing that
$\phi(x)\in
\bnp$. Hence $\phi(\bnp)\subset \bnp$.
 \end{proof}

For some groups, such as free groups and free-solvable groups
$F/F^{\sss (n)}$, the derived series and the torsion-free derived
series coincide.

\begin{proposition}\label{torsionfree} If $G$ is a group such that,
for each $n$,
$G^{\sss (n)}/G^{\sss (n+1)}$ is torsion-free as a $\Z[G/G^{\sss
(n)}]$-module, 
 then the torsion-free derived series of $G$ agrees with
the derived series of $G$. Hence for a free group $F$, $F^{\sss
(n)}_{\sss H}=F^{\sss (n)}$ for each $n$.
\end{proposition}

\begin{proof} By definition,
$\gn=G^{\sss (0)}=G$. Suppose
$\gn=G^{\sss (n)}$. Then, under
the hypotheses, $\gn/[\gn,\gn]$ is a torsion-free module and hence
$\gnp = \ker\pi_n=[\gn,\gn]=[G^{\sss (n)},G^{\sss (n)}]=G^{\sss
(n+1)}$.

It is well known that $F^{\sss (n)}/F^{\sss (n+1)}$ is a $\Z
[F/F^{\sss (n)}]$-torsion-free
module. This can be seen by examining the free $\Z[F/F^{\sss (n)}]$
cellular chain complex for the covering space of a wedge of circles
corresponding to the subgroup
 $F^{\sss (n)}$. The module $F^{\sss (n)}/F^{\sss (n+1)}$ is merely
the first homology of
this chain complex. Since the chain complex can be chosen to have
no $2$-cells, its first homology
 is a submodule of a free module and thus is torsion-free. Hence, by
the first part of this proposition, the derived
series and the
 torsion-free derived series of a free group agree.
 \end{proof}

In \cite{CH}, T. Cochran and the author prove a version of Stallings' Theorem \cite{St} for the derived series. Specifically, it was there shown that if 
a map of finitely presented groups $\phi: A \rightarrow B$ is rationally $2$-connected then it induces a monomorphism 
$$\phi_{\ast}: A/A_{\sss H}^{\sss (n)} \hookrightarrow B/B_{\sss H}^{\sss (n)}$$ for all $n\geq 0$ \cite[Theorem 4.1]{CH} (see Theorem~\ref{CH_thm} and Proposition~\ref{CHprop} below).
For this paper, we will need the following ``generalization'' of that theorem. The following theorem is in fact a consequence of the proof of Theorem 4.1 of \cite{CH}. 
For the convenience of the reader, we will sketch the proof of Theorem~\ref{CH_thm} after the statement of Proposition~\ref{CHprop}.

\begin{theorem}[Scholium to Theorem 4.1 of \cite{CH}]\label{CH_thm}
Let $n$ be a non-negative integer or $n=\omega$.  If $\phi: A \rightarrow B$ is a homomorphism that induces a monomorphism
$\phi_\ast: H_1(A;\mathcal{K}(B/\bk))\to H_1(B;\mathcal{K}(B/\bk))$ for each $k\leq n$, then for each
$k\leq n$, $\phi$ induces a monomorphism $A/\akp
\hookrightarrow
B/\bkp$.  
Moreover, if $\phi: A \rightarrow B$ induces an isomorphism
$\phi_\ast: H_1(A;\mathcal{K}(B/\bk))\to H_1(B;\mathcal{K}(B/\bk))$ for each $k\leq n$, then for each
$k\leq n$, $\phi$ induces a monomorphism
$\ak/\akp
\hookrightarrow
\bk/\bkp$ between modules of the same rank (over $\Z[A/\ak]$ and $\Z[B/\bk]$ respectively).  In addition, if $\phi$ is onto then 
$\phi_{\ast}: A/\akp
\rightarrow
B/\bkp$ is an isomorphism.
\end{theorem}

The following proposition guarantees that one of the hypotheses of Theorem~\ref{CH_thm} is satisfied whenever $\phi$ is a rationally $2$-connected map. 
Note that Proposition~\ref{CHprop} and Theorem~\ref{CH_thm} together imply Theorem~4.1 of \cite{CH}. The justification for calling Theorem~\ref{CH_thm} a generalization of 
Theorem~4.1 of \cite{CH},
is that, in the subsequent sections we will describe several conditions (see Proposition~\ref{prop:amal_prod}, Theorem~\ref{diff_rho}, Theorem~\ref{rho0}, Proposition~\ref{rhoadd}, and Lemma~\ref{retr})
under which the hypothesis of Theorem~\ref{CH_thm} is satisfied but where the $2$-connected hypothesis of Theorem~4.1 of \cite{CH} (and Proposition~\ref{CHprop} and
Proposition~4.3 of \cite{CH}) fails.

\begin{proposition}[Proposition 4.3 of \cite{CH}] \label{CHprop} Let $A$ be a
finitely-generated group and $B$ a finitely related group. Suppose
$\phi:A\to B$ induces
a monomorphism (resp. isomorphism) on $H_1(- ;\mathbb{Q})$ and an
epimorphism on $H_2(- ;\mathbb{Q})$. Then for each $k\geq 0$,  $\phi$
induces a monomorphism (resp. isomorphism)
$\phi_\ast: H_1(A;\mathcal{K}(B/\bk))\to
H_1(B;\mathcal{K}(B/\bk))$.
\end{proposition}

\begin{proof}[Proof of Theorem~\ref{CH_thm}] We sketch the inductive proof of the first claim of the theorem, referring the reader to \cite{CH} for more details. For $n=0$,
 $A/A^{\sss (1)}_{\sss H}$ is merely
$H_1(A;\mathbb{Z})/\{\mathbb{Z}\text{-Torsion}\}$. But $A/A^{\sss (0)}_{\sss H}=\{e\}$ so $\mathcal{K}(A/A^{\sss (0)}_{\sss H})=\mathbb{Q}$. Thus our hypothesis, that $\phi$ induces
a monomorphism on $H_1(-;\mathbb{Q})$, implies that $\phi$
induces a monomorphism on $H_1(-;\mathbb{Z})$ modulo torsion.
Now assume that $\phi$ induces a monomorphism $A/\an \subset
B/\bn$. We will prove that this
holds for $n+1$.

It follows from Proposition~\ref{functoriality} that $\phi(\anp)\subset \bnp$. Thus from the commutative diagram below we see that it suffices to show that
$\phi$ induces a monomorphism $\an/\anp\to\bn/\bnp$.
\begin{diagram}
1      &\rTo &    \an/\anp  &\rTo &  A/\anp  &\rTo  & A/\an   & \rTo &1\\
 & & \dTo_{\phi}  &   & \dTo_{\phi_{n+1}} & & \dTo_{\phi_{n}} & &\\
 1      &\rTo&    \bn/\bnp  &\rTo&    B/\bnp  &\rTo & B/\bn    & \rTo & 1
\end{diagram}
Now suppose that $\an/\anp\to\bn/\bnp$ were not injective. From our discussions above in the proof of Proposition~\ref{functoriality} we see that there would exist an $a\in\an$ representing a non-torsion class $[a]$ in ${\an}/{[\an,\an]}$ such that $\phi (a)$
represents a torsion class in ${\bn}/{[\bn,\bn]}$. But
$$
\an/[\an,\an]\cong H_1(A;\mathbb{Z}[A/\an]).
$$
The torsion submodule is characterized precisely as the kernel of the canonical map
$$
H_1(A;\mathbb{Z}[{A}/{\an}])\to H_1(A;\mathbb{Z}[A/\an])\otimes_{\mathbb{Z}[A/\an]}\mathcal{K} (A/\an)\cong H_1(A;\mathcal{K}(A/\an)).
$$
A similar statement holds for $B$. But the inductive hypothesis that $A/\an \subset B/\bn$ guarantees that
$$
H_1(A;\mathcal{K}(A/\an))\to H_1(A;\mathcal{K}(B/\bn))
$$
is injective.  Moreover, the hypothesis of the theorem guarantees that
$$
H_1(A;\mathcal{K}(B/\bn))\to H_1(B;\mathcal{K}(B/\bn))
$$
is injective, leading to a contradiction.
\end {proof}

\section{Definition of $\rhon$}\label{section:def_rhon}

On the class of closed, oriented $(4k-1)$-dimensional manifolds, we will define a $\Q$-homology cobordism invariant $\rho_{n}$ for each $n\in \mathbb{N}$.   This will be defined as a signature defect associated to 
the $n^{th}$ term of the torsion-free derived series of the fundamental group of the manifold.  We begin by 
recalling the definition of the $L^{2}$-signature of a
$4k$-dimensional manifold.  For more information on $L^{2}$-signature and $\rho$-invariants see \cite[Section 2]{CT}, \cite[Section 5]{COT} and  \cite{LS}.

Let $\Lambda$ be a countable group and $\mathcal{U}\Lambda$ be the algebra of unbounded operators affiliated to $\mathcal{N}\Lambda$, the von Neumann algebra of $\Lambda$.
Then $\sigma^{(2)}_{\Lambda}: \text{Herm}_{n}(\mathcal{U}\Lambda) \rightarrow \R$ is defined by 
$$\st_{\Lambda}(h)=\text{tr}_{\sss
\Lambda}(p_+(h))-\text{tr}_{\sss
\Lambda}(p_-(h))$$ 
for any  $h \in \text{Herm}_{n}(\mathcal{U}\Lambda)$ where
$\text{tr}_{_\Lambda}$ is the von Neumann trace and $p_\pm$ are the
characteristic functions on the positive and negative reals.  It is known that $\st_{\Lambda}$ can be extended to the Witt group of Hermitian forms on finitely generated projective
$\mathcal{U}\Lambda$-modules.  

\begin{lemma}[see for example Corollary~5.7 of \cite{COT} and surrounding discussion]\label{witt} 
The $L^{2}$-signature, $\sigma^{(2)}$,
 is a well-defined real-valued homomorphism on the Witt group of hermitian forms on finitely 
generated projective $\mathcal{U}\Lambda$-modules. Restricting this homomorphism to nonsingular forms on free modules gives 
$$\st_{\Lambda}: L^{0}(\mathcal{U}\Lambda) \rightarrow \R.$$
\end{lemma}

\noindent In particular, if $h$ is a nonsingular pairing with a metabolizer then $\st_{\Lambda}(h)=0$.

Let $W$ be $4k$-dimensional manifold and $\Phi: \pi_1(W) \rightarrow
\Lambda$ be a coefficient system for $W$.  Let $h_{W,\Lambda}$ be the
composition
of the following homomorphisms 
\begin{equation}\label{eq:herm}H_{2k}(W;\mathcal{U}\Lambda) \rightarrow H_{2k}(W,\partial
W;\mathcal{U}\Lambda) \xrightarrow{\text{PD}} 
H^{2k}(W;\mathcal{U}\Lambda) \xrightarrow{\kappa}
H_{2k}(W;\mathcal{U}\Lambda)^\ast
\end{equation}
where $H_{2k}(W;\mathcal{U}\Lambda)^\ast =
\text{Hom}_{\mathcal{U}\Lambda}(H_{2k}(W;\mathcal{U}\Lambda),\mathcal{U}\Lambda)$.  Since $\mathcal{U}\Lambda$ is a von Neumann regular ring,
 the modules $H_{2k}(W;\mathcal{U}\Lambda)$ are finitely generated 
projective right $\mathcal{U}\Lambda$-modules.  Then $h_{W,\Lambda} \in \text{Herm}_{n}(\mathcal{U}\Lambda)$
and we define $\st(W,\Lambda)=\st_{\Lambda}(h_{W,\Lambda})$.
We will sometimes write $\st(W,\Lambda)$ as $\st_\Lambda(W)$ or
$\st(W,\Phi)$ when we want to emphasize the map $\Phi$.
 
Suppose that $\Lambda$ is PTFA. Let $U$ be a (possibly empty) union of components of the boundary of $W$.  Then $\Z\Lambda$ embeds in its right ring of quotients $\mathcal{K}\Lambda$.  Moreover, the
map from $\Z\Lambda$ to $\mathcal{U}\Lambda$ factors as $\Z\Lambda \rightarrow \mathcal{K}\Lambda \rightarrow \mathcal{U}\Lambda$ 
making $\mathcal{U}\Lambda$ into a $\mathcal{K}\Lambda-\mathcal{U}\Lambda$-bi-module.  Since any module over a skew field is free, $\mathcal{U}\Lambda$ is a flat $\mathcal{K}\Lambda$-module.
Hence, $H_{2k}(W,U;\mathcal{U}\Lambda)\cong H_{2}(W,U;\mathcal{K}\Lambda)\otimes_{\mathcal{K}\Lambda}\mathcal{U}\Lambda$. In particular, $H_{2k}(W,U;\mathcal{K}\Lambda)=0$ if and only if
$H_{2}(W,U;\mathcal{U}\Lambda)$=0.

We will use the following facts about $L^{2}$-signatures throughout this paper.  The first two remarks follow directly from the definition of $\sigma^{(2)}(W,\Lambda)$.

\begin{remark}  Suppose $W$ is a compact, oriented $4k$-dimensional manifold and $\Psi: \pi_{1}(W)\rightarrow \Lambda$ is a coefficient system for
$W$. \label{fiveremarks}
\begin{enumerate}

\item \label{remark:from_bound} If $H_{2k}(W,U;\mathcal{U}\Lambda)=0$ (or $H_{2k}(W,U;\mathcal{K}\Lambda)=0$ if $\Lambda$ is PTFA), 
where $U$ is a (possibly empty) union of components of the boundary of $W$, then $\sigma^{(2)}(W,\Lambda)=0$. 

\item \label{remark:EasyCase}If $\Lambda=\{1\}$ is the trivial group then $\sigma^{(2)}(W,\Lambda)=\sigma(W)$ where $\sigma$ is the ordinary signature function.

\item 
\label{gamma_induction} If $\Lambda \subset \Lambda^\prime$ then
$\st(W,\Lambda)=\st(W,\Lambda^{\prime}) $ (see for example, Proposition~5.13 of \cite{COT}).

\item \label{addonbound} 
Suppose $V$ is another compact, oriented $4k$-dimensional manifold, $\Psi^{\prime}:\pi_{1}(V)\rightarrow \Lambda$ is a coefficient system for $V$ and $(V,\Psi^{\prime})$ has the same oriented boundary as $(W,\Psi)$ (meaning the maps to $\Lambda$ agree on the boundary)
then
\[\st(W \cup_{\partial W}\overline{V}, \Psi \cup \Psi^{\prime})
= \st(W,\Psi)- \st(V,\Psi^{\prime})\](see for example, Lemma~5.9 of \cite{COT}).

\item \label{signature_zero} If 
$W$ is closed then $\st(W,\Lambda)=\sigma(W)$ (see for example, Lemma~5.9 of \cite{COT}).
\end{enumerate}
\end{remark}

We now define $\rho_\Gamma(M)$ for a $(4k-1)$-dimensional manifold and coefficient system 
$\pi_{1}(M)\rightarrow \Gamma$.  Let $M$ be a closed, orientable, $l$-dimensional manifold with $l
\not\equiv 0 \mod 4$.  It is well known that $rM$ is the
boundary of some compact, orientable manifold for $r\in \{1,2\}$.
J.-C. Hausmann showed further \cite[Theorem~5.1]{Haus} that
$rM$ is the boundary of a compact, orientable manifold $W$ for which
the inclusion map of $M$ into 
$W$ induces a monomorphism on $\pi_1$.  That is, $rM$ is \emph{stably
nullbordant over $\pi_1(M)$} for some $r\in \{1,2\}$ in the language of Definition~\ref{def:stably} below.  In \cite{CW}, S.
Chang and S. Weinberger use this fact to define 
a new ``Hirzebruch type'' invariant $\tau_{(2)}$ for a $(4k-1)$
dimensional manifold $M$ by setting
$\tau_{\sss (2)}(M)= \frac{1}{r}(\st(W,\pi_1(W))-\sigma(W))$.  To
define $\rho_n$ we proceed in a similar manner, stably bounding
over $\pi_1(M)_n$ instead of $\pi_1(M)$.  To do this, we will show
that $rM$ is 
stably nullbordant over $\Gamma$ for any coefficient system
$\pi_1(M) \rightarrow \Gamma$.  

\begin{definition}\label{def:stably}Let $M=M_1 \cup \dots \cup M_m$ be a disjoint union
of $m$ connected, closed, oriented $l$-dimensional manifolds and
$\mathcal{S} = \{\phi_i:\pi_1(M_i) \rightarrow \Gamma_i\}_{i=1}^m$
be a collection of coefficient systems for $M$.   
We say that $(M,\mathcal{S})$ is \textbf{stably nullbordant} (or
\textbf{s-nullbordant}) if  there exists a triple
$(W,\Phi,\mathcal{T})$ 
where $W$ is a compact, connected, oriented $(l+1)$-dimensional
manifold with $\partial W=M$,  
$\Phi: \pi_1(W) \rightarrow \Lambda$ is a coefficient system for $W$,
and $\mathcal{T}=\{\theta_i:\Gamma_i \into \Lambda\}_{i=1}^m$
 is a collection of monomorphisms such that for each $1\leq i\leq m$,
the following diagram commutes (after modifying $\Phi$ by a change of basepoint isomorphism)
\begin{equation}\label{diag1}
\begin{diagram} \pi_1(M_i) & \rTo^{\phi_i} & \Gamma_i \\
\dTo_{(\iota_i)_\ast} && \dIntoA_{\theta_i} \\
\pi_1(W) & \rTo^\Phi & \Lambda
\end{diagram}
\end{equation}
\noindent where $\iota_i: M_i \rightarrow W$ is the inclusion map.
We call the triple $(W,\Phi,\mathcal{T})$ a \textbf{stable (or s-)nullbordism}
for $(M,\mathcal{S})$.
We say that $(M_1,\mathcal{S}_1)$ is \textbf{stably (or s-)bordant} to
$(M_2,\mathcal{S}_2)$ if $(M_1 \cup \overline{M}_2, \mathcal{S}_1
\cup \mathcal{S}_2)$ is s-nullbordant.
If, in addition, $\Gamma_i \cong \Gamma$ for each $i=1,\dots,m$, we
say that $M$ is \textbf{stably (or s-)nullbordant over $\Gamma$} or that
$(M,\Gamma)$ is \textbf{stably (or s-)nullbordant}.
\end{definition}

We remark that stable bordism is an equivalence
relation since if $(W_i,\Phi_i: \pi_1(W_i)\rightarrow \Lambda_i,
\{\theta^i_i,\theta^i_{i+1} \})$ is an s-nullbordism
for $(M_i\cup \overline{M}_{i+1},\{\phi_i,\phi_{i+1}\})$ ($i=1,2$)
then $(W=W_1 \cup_{M_2} W_2,\Phi:\pi_1(W) \rightarrow \Lambda_1
\ast_{\Gamma_2} \Lambda_2,\{\theta^\prime_1,\theta^\prime_3\})$ is an
s-nullbordism for $(M_1\cup \overline{M}_3,\{\phi_1,\phi_3 \} )$
when $\Phi=\Phi_1 \ast \Phi_2$, and $\theta_{1}^{\prime}$, $\theta_{3}^{\prime}$
are the obvious compositions.

The proof of the following lemma is similar to the proof of Hausmann's Theorem~5.1.
\cite{Haus}.

\begin{lemma}\label{lem:sbord}Let $M=M_1 \cup \cdots \cup M_m$ be a disjoint union of
closed, connected, oriented $l$-dimensional manifolds and 
$\mathcal{S} = \{\phi_i:\pi_1(M_i) \rightarrow \Gamma_i\}_{i=1}^m$ be
a collection of coefficient systems.   
If $M$ is nullbordant then $(M,\mathcal{S})$ is s-nullbordant.\end{lemma}
\begin{proof}If $\Gamma = \Gamma_1 * \cdots * \Gamma_m$ is the free
product of the collection $\{\Gamma_i\}$ then there is 
a natural inclusion $\Gamma_i \into \Gamma$ for each $i$.   
By the homological coning construction of W. Thurston-W. Kan,
$\Gamma$ is subgroup of an acyclic group $\Lambda_\Gamma$ \cite[Section 3]{KT}.  For each
$i$, let $\theta_i: \Gamma_i \into \Lambda_\Gamma$ be the inclusion 
of $\Gamma_i$ into the acyclic group $\Lambda_\Gamma$.  The
collection $\{\theta_i \circ \phi_i\}$ gives us a map $f:M
\rightarrow K(\Lambda_\Gamma,1)$ such that 
$(f_{| M_i})_\ast = \theta_i \circ \phi_i$ for each $i$ and 
the following diagram commutes.
\[\begin{diagram} M & \rTo^{f\hspace{10pt}} & K(\Lambda_\Gamma,1) \\
\dTo_{\iota} & & \dTo \\
W & \rTo & \text{pt}
\end{diagram}
\]
Since $K(\Lambda_\Gamma,1)$ is acyclic, the map $K(\Lambda_\Gamma,1)
\rightarrow \text{pt}$ induces an isomorphism on integral
homology, hence an isomorphism on oriented bordism theory.  
Since $M =\partial W$,
this implies that there is a compact manifold $W^\prime$ 
and map $g:W^\prime \rightarrow K(\Lambda_\Gamma,1)$ such that
$g_{|\partial W^{\prime}}=f$.  Hence $(W^\prime,g_\ast, \{\theta_i\})$ is an
s-nullbordism for $M$.  
\end{proof}

Since the $l$-dimensional oriented bordism group
$\Omega_l^\text{or}(\text{pt})$ is 2-torsion when $l \not\equiv
0 \mod 4$, there is always an $(l+1)$-dimensional manifold $W$ such that 
$\partial W=2M$.  Hence we see that $(2M,\{\phi\})$ is always
s-nullbordant.

\begin{corollary}\label{cor:sbord}If $M$ is a closed, connected, oriented
$l$-dimensional manifold with $l \not\equiv 0 \mod 4$ and $\phi:\pi_1(M)\rightarrow \Gamma$ is a coefficient system
for $M$ then there is 
an integer $r\in \{1,2\}$ such $rM$ is s-nullbordant over $\Gamma$. 
\end{corollary}

\noindent In this paper, we will often assume that $M$ is a $3$-dimensional
manifold.  Since every closed, oriented $3$-dimensional manifold is the boundary of a
compact, oriented $4$-dimensional manifold, $(M,\phi)$ is s-nullbordant
for any $\phi: \pi_1(M)\rightarrow \Gamma$. 

For a $(4k-1)$-dimensional closed, oriented manifold $M$ and homomorphism $\phi: \pi_{1}(M)\rightarrow \Gamma$,
we define $$\rhot(M,\phi):= \frac{1}{r}(\st(W,\Lambda) - \sigma(W))$$ 
for $(W,\Psi)$ any s-nullbordism for $(rM,\phi)$ where $\sigma(W)$ is the ordinary signature of $W$.  
By the following Lemma, this definition only depends on $M$ and $\phi$.

\begin{lemma}\label{lem:independence}$\rhot(M,\phi)$ is independent of the choice of
$(W,\Phi,\theta)$. 
\end{lemma}
\begin{proof}Let $(W,\mathbf{\Phi},\mathcal{T})$ and
$(W^\prime,\Phi^\prime,\mathcal{T}^\prime)$ be two s-nullbordisms for $(r_{1}M,\{\phi\})$ and  $(r_{2}M,\{\phi\})$
respectively.  Assume that $r_{1}=r_{2}=1$.
Let $C=W\cup \overline{W^\prime}$ be the closed, oriented
$4k$-manifold obtained by gluing $W$ and $W^\prime$ along $M$.  Then 
there is a coefficient system for C, 
$$\Phi \ast \Phi^\prime:
\pi_1(C)=\pi_1(W)\ast_{\pi_1(M)}\pi_1(W^\prime) \rightarrow \Lambda
\ast_{\pi_1(M)}\Lambda^\prime =\Lambda \ast_{\Gamma} \Lambda^{\prime}$$
such that $\Phi\ast\Phi^\prime(\alpha)=i_\Lambda(\Phi(\alpha))$ for
all $\alpha \in \pi_1(W)$ where   
$i_\Lambda: \Lambda \rightarrow \Lambda \ast_\Gamma \Lambda^\prime$
sends $\lambda \in \Lambda$ to the word $\lambda \in \Lambda
\ast_\Gamma \Lambda^\prime$ (similarly for $\alpha \in
\pi_1(W^\prime)$).
Since $\theta: \Gamma \rightarrow \Lambda$ and $\theta^\prime:\Gamma
\rightarrow \Lambda^\prime$ are monomorphisms, the maps
$i_\Lambda$ and $i_{\Lambda^\prime}$ are monomorphisms.  Hence by
Remark~\ref{fiveremarks}~(\ref{gamma_induction}), 
$\st(W,\Lambda)=\st(W,\Lambda\ast_\Gamma \Lambda^\prime)$ and
$\st(W^\prime,\Lambda^\prime)=\st(W,\Lambda\ast_\Gamma
\Lambda^\prime).$
Moreover, by Remark~\ref{fiveremarks}~(\ref{addonbound})
\begin{equation}\label{uselhs}\st_{\Lambda \ast_\Gamma
\Lambda^\prime}(C)-\sigma(C)=(\st_{\Lambda \ast_\Gamma
\Lambda^\prime}(W)-\sigma(W))-(\st_{\Lambda \ast_\Gamma
\Lambda^\prime}(W^\prime)-\sigma(W^\prime)).\end{equation}   
Since $C$ is closed, by Remark~\ref{fiveremarks}~(\ref{signature_zero}), the left hand side of (\ref{uselhs}) is
0.  Therefore,
$\st_\Lambda(W)-\sigma(W)=\st_{\Lambda^\prime}(W^\prime)-\sigma(W^\prime)$.  The proofs when $r_{i}\neq 1$
are similar and are not included since most of our applications focus on $3$-manifolds where $r$ can be assumed to be $1$.\end{proof}

Note that we may occasionally write $\rhot(M,\phi)$ as
$\rhot_\Gamma(M)$ or $\rhot(M,\Gamma)$ when the map $\phi$ is clear. 
As a result of Remark~\ref{fiveremarks}~(\ref{gamma_induction}), $\rhot_\Gamma(M)$ only depends on the image of
$\pi_1(M)$ in $\Gamma$.

\begin{lemma}[$\Gamma$-induction]\label{gamma_induction_for_rho}Suppose
$M$ is a closed, oriented $(4k-1)$-dimensional manifold and $\phi:
\pi_1(M) \rightarrow \Gamma$ is a
coefficient system for $M$. If $\iota: \Gamma \hookrightarrow
\Gamma^\prime$ is a monomorphism then
$\rhot(M,\phi)=\rhot(M,\iota \circ \phi)$.
\end{lemma}

\begin{definition}For each $0\leq n\leq \omega$ and $(4k-1)$-dimensional closed, oriented manifold $M$,  
we define the \textbf{$n^{th}$-order $\rho$-invariant of $M$} by $$\rho_n(M):=\rho(M, \phi_n:G \onto G/G_{\sss H}^{\sss (n+1)}) \in \R$$
where $G=\pi_{1}(M)$.\end{definition}

We will now define the $n^{th}$-order $\rho$-invariant of a link in $S^{3}$.  First, suppose $L \subset S^3$ is an 
$m$-component link in $S^3$ with linking numbers $0$.  Let $\text{N}(L)$ be a neighborhood of $L$ in $S^3$.

 \begin{proposition}\label{linkgroup}Let $L \subset S^{3}$ be a link for which all the pairwise linking numbers are zero and $G=\pi_{1}(S^{3} - \text{N}(L))$.  
 The  longitudes of $L$ lie 
in $G^{\sss (\omega)}_{\sss H}$.
\end{proposition}
\begin{proof} Let  $\lambda_i$ be the longitude of the $i^{th}$ component of $L$. We will show that for each $n\geq 1$, $\lambda_{i} \in \gn$.
Since the linking numbers of $L$ are zero, $\lambda_{i} \in \gone$. 
Suppose that $\lambda_i\in
\gn$, for some $n\geq 1$. Since the longitudes lie 
on the boundary tori, they commute with the meridians $x_i$. Hence for each $i$, $[x_i,\lambda_i]=1$ is a relation in $G$.
The relation $[x_i,\lambda_i]=1$ in $G$ creates the
relation $\lambda_{i}(1-x_i)=0$ in the module $\gn /[\gn ,\gn]$,
showing that $\lambda_i\in \gnp$ since $x_{i}\neq 1$ in $G/\gn$. 
\end{proof}

By contrast, 
the longitudes rarely lie in $G_{\omega}$, much less lie in
$G^{\sss (\omega)}$, the former being true if and only if all 
of Milnor's $\overline{\mu}$-invariants are zero. The Borromean
Rings and the Whitehead links provide examples 
where the longitudes lie in $G^{\sss (\omega)}_H$ but not in $G_{\omega}$.

As a corollary, performing $0$-framed surgery on a link with linking numbers $0$ does not change the quotient of the link group by
a term its torsion-free derived series.

\begin{corollary}[see Proposition 2.5 of \cite{CH}] Suppose $L\subset S^{3}$ is a link with linking numbers $0$.  Let $M_L$ be the closed $3$-manifold obtained by
performing $0$-framed surgery on the components of $L$ and 
$i:S^3 - \text{N}(L) \rightarrow M_L$ be the inclusion map.   Then for
each $n \geq 1$, 
\[\frac{\pi_1(S^3 - \text{N}(L))}{\pi_1(S^3 - \text{N}(L))_H^{(n)}}
\xrightarrow{\cong} 
\frac{\pi_1(M_L)}{\pi_1(M_L)_H^{(n)}} 
\]
\end{corollary}
\begin{proof}The kernel of $i_{\ast}$ is the normal subgroup generated by the longitudes.   But by Proposition~\ref{linkgroup} above, the longitudes lie in 
$\pi_{1}(S^{3} - \text{N}(L))^{\sss (n)}_{\sss H}$ for all $n\geq 0$.
\end{proof}

\noindent  Hence it makes sense to make the following definition 
of $\rho_{n}$ for a link.  Note that the following definition does not require that the linking numbers be $0$.

\begin{definition}Let $L$ be a link in $S^3$.  For each $0 \leq n \leq \omega$ we
define \[\rho_n(L)=\rho \left(M_L, \phi_n:\pi_1(M_L) \onto
\frac{\pi_1(M_L)}{\pi_1(M_L)_H^{(n+1)}} \right).\]
\end{definition}

\noindent As an easy example, we show that $\rho_{n}(\#_{i=1}^m S^1 \times S^2)=\rho_{n}(\text{trivial link})=0$
for all $0 \leq n \leq \omega$.

\begin{example}\label{example:TrivLink}
Let $W$ be the boundary connected sum of $m$ copies of $S^1 \times
D^3$. Then $\partial W=\#_{i=1}^m S^1 \times S^2$. 
Moreover, the inclusion  $i:\partial W \rightarrow W$ induces an
isomorphism on $\pi_1$.
 Let $\Gamma_n = F/F_H^{(n+1)}$ where $F=\pi_1(W)$ is the free group
with $m$ generators. By definition, $\rho_n(\partial
W)=\st_{\Gamma_n}(W)-\sigma(W)$.
Since $W$ is homotopy equivalent to a 1-complex,
$H_2(W;\Z\Gamma_n)=0$.  
Hence, $\st_{\Gamma_n}(W)=\sigma(W)=0$. In particular, 
$\rho_{n}(\#_{i=1}^m S^1 \times S^2)=\rho_{n}(\text{trivial link})=0$.
\end{example}

The most important and easiest example to understand is $\rho_{0}$ for a knot in $S^{3}$.  In this case, $\rho_{0}$ is determined by the Levine-Tristram signatures of the knot.

\begin{example}\label{ex:knot}Let $K$ be a knot in $S^{3}$.  By Lemma~5.4 of \cite{COT} and Lemma~5.3 of \cite{COT1}, $$\rho_{0}(K)=\int_{S^{1}} \sigma_{\omega}(K)\text{d}\omega$$
where  $\sigma_{\omega}(K)$ is the
Levine-Tristram signature of $K$ at $\omega\in S^{1}$ and the circle is normalized to have length $1$.  Since $\beta_{1}(M_{K})=1$,  the Alexander module of $M_{K}$ is torsion \cite[Proposition 2.11]{COT}.  
Therefore, for $n\geq 0$, $\pi_{1}(M_{K})^{\sss (n+1)}_{\sss H}=\pi_{1}(M_{K})^{(1)}$; hence $\rho_{n}(K)=\rho_{0}(K)$.
\end{example}

A nice property of $\rho_n$ is that it is additive under the
connected sum of manifolds.  To prove this, we show
that the torsion free derived series behaves well under the
inclusion $A \rightarrow A \ast_{\sss C} B$ for suitable $C$.
We start with a lemma.

\begin{lemma}\label{lemma:amal_prod}Let $\phi :C \rightarrow G$ be a
homomorphism. 
 If $\beta_1(C)=0$ then the image of $\phi$ is contained in $G_{\sss
H}^{\sss (\omega)}$.
\end{lemma}
\begin{proof}We will show by induction on $n$ that $\phi(C) \subset
\gn$ for all $n\geq 0$.  This is trivial when $n=0$.  
Suppose for some $n\geq 0$, that $\phi(C) \subset \gn$.  Then $\phi$
induces a map 
$\phi_\ast: C/[C,C] \rightarrow \gn/[\gn,\gn] \rightarrow \gn/\gnp$.
Since $\beta_1(C)=1$, $C/[C,C]$ is $\Z$-torsion.  However, 
$\gn/\gnp$ is $\Z$-torsion free, hence $\phi_\ast$ is trivial which
implies $\phi(C) \subset \gnp$.
\end{proof}

\begin{proposition}\label{prop:amal_prod}Let $A\ast_{\sss C} B$ be
the amalgamated product of $A$ and $B$ where 
$C \into A$ and $C \into B$ are monomorphisms and $\beta_1(C)=0$.
For each $0 \leq n\leq \omega$, the inclusion $i : A \rightarrow A
\ast_{\sss C} B$ induces  
a monomorphism 
\begin{equation}\label{eq:amal_prod}
 i_\ast : \frac{A}{A^{(n)}_H} \into \frac{A \ast_{\sss C}
B}{(A\ast_{\sss C} B)_H^{(n)}}.\end{equation}
\end{proposition}
\begin{proof}Let $G=A*_{\sss C} B$.  For each $n\geq 0$, we have the
following Mayer-Vietoris sequence for group homology with
$\mathcal{K}G_n$-coefficients 
\[ \rightarrow H_1(C;\mathcal{K}G_n) \rightarrow H_1(A;\mathcal{K}G_n) \oplus
H_1(B;\mathcal{K}G_n) \rightarrow H_1(G;\mathcal{K}G_n)
\rightarrow  \]
where the coefficients systems for $A,B,C$ and $G$ are the obvious
ones. By Lemma~\ref{lemma:amal_prod}, the image
of 
$C \rightarrow G \onto G_n$ is trivial hence $H_1(C;\Z G_n)\cong
H_1(C;\Z)\otimes_\Z \Z G_n$.  

Since $\Q$ is a flat $\Z$-module, $\mathcal{K}G_n$ is a 
flat $\Z G_n$-module, and $\beta_1(C)=0$, it follows that
$H_1(C;\mathcal{K}G_n)\cong H_1(C;\Q) \otimes_\Q \mathcal{K}G_n =
0$. Here, the map $\Q \rightarrow \mathcal{K}G_n$ is
induced by $1 \into G_n$.  Thus 
$i_\ast : H_1(A;\mathcal{K}G_n) \rightarrow
H_1(G;\mathcal{K}G_n)$ is a monomorphism.    By
Theorem~\ref{CH_thm}, $i$ induces a monomorphism $A/\an \into
G/\gn$ as desired.
Since $A/\an \into G/\gn$ for all $n\geq 0$, it follows immediately
from the definition of $A_{\sss H}^{\sss (\omega)}$ and $G_{\sss
H}^{\sss (\omega)}$ that 
$A/A_{\sss H}^{\sss (\omega)} \into G/G_{\sss H}^{\sss (\omega)}$. 
\end{proof}

\begin{proposition}\label{prop:connsum}Let $k\geq 1$ and let $M_1$ and $M_2$ be closed, oriented, connected
$(4k-1)$-dimensional manifolds.  For each $0 \leq n\leq \omega$, 
\[\rhon(M_1 \# M_2) = \rhon(M_1) + \rhon(M_2).\]
\end{proposition}
\begin{proof} Let $W$ be the $4k$-manifold obtained by adding a 1-handle
to $(M_1 \sqcup M_2)\times I$ along some $D^{4k-1} \sqcup D^{4k-1} \subset (M_1
\sqcup M_2)\times \{1\}$ and $G=\pi_1(W)$.  Then
$\partial W = M_1 \sqcup M_2 \sqcup \overline{M_1 \# M_2}$ so 
 $\st(W,G_n)-\sigma(W) = \rho(M_1,G_n) + \rho(M_2, G_n) -
\rho(M_1 \# M_2, G_n)$.
Since the inclusion $i: \overline{M_1 \# M_2} \rightarrow W$ induces
an isomorphism on $\pi_1$, by Lemma~\ref{gamma_induction_for_rho},
$\rhon(\overline{M_1 \# M_2})=\rho(\overline{M_1 \# M_2}, G_n)$.
Moreover $G = \pi_1(M_1) \ast \pi_1(M_2)$ and the inclusion
map $i_1:M_1 \rightarrow W$ induces the inclusion map $(i_1)_\ast :
\pi_1(M_1) \rightarrow \pi_1(M_1) \ast \pi_1(M_2)$ on $\pi_1$.   
Therefore, by Proposition~\ref{prop:amal_prod}, $(i_1)_\ast :
\pi_1(M_1)_n \to G_n$ is a monomorphism for all $n\leq \omega$.
Thus, by Lemma~\ref{gamma_induction_for_rho},
$\rhon(M_1)=\rho(M_1,G_n)$ (similarly for $M_2$).   

To finish the proof, it suffices to show that
$\st(W,G_n)-\sigma(W)=0$.  To see this, we point out that $(W,M_1 \#
M_2)$ is 
homotopy equivalent to $(W^\prime,M_1 \# M_2)$ where $W^\prime$ is
obtained by attaching a $(4k-1)$-dimensional cell to $M_1 \# M_2$.  Therefore,
$H_{2k}(W,M_1 \# M_2;\Lambda)=0$ for any coefficient system $\phi:
\pi_1(W) \rightarrow \Gamma$ and left $\Z\Gamma$-module $\Lambda$.   Therefore, $H_{2k}(W,M_1 \# M_2;\mathcal{U}G_{n})=0$
and hence
$\st(W,G_n)=\sigma(W)=0$.
\end{proof}

\section{Homology Cobordism}\label{section:homcob}

Our primary interest in this paper is the study of manifolds up to (rational) \emph{homology cobordism}.  We begin with a definition. 

\begin{definition}Let $M^m_1$ and $M^m_2$ be oriented, closed
$m$-dimensional manifolds.  We say that $M_1$ is \textbf{$\Q$-homology cobordant} (respectively \textbf{homology cobordant})
to $M_2$ if there exists a  oriented, $(m+1)$-manifold
$W^{m+1}$ such that $\partial W = M_1 \cup \overline{M_2}$, and the 
inclusion maps $i_j : M_j \rightarrow W$ induce isomorphisms on
$H_{\ast}(-;\Q)$ (respectively $H_{\ast}(-;\Z)$).  In this case we write $M_1 \sim_{\Q H} M_2$ (respectively $M_1 \sim_{\Z H} M_2$) and 
define the set of rational (respectively integral) homology cobordism classes of $m$-dimensional manifolds
to be $\mathcal{H}^m_\Q = \{M^m\}/\sim_{\Q H}$ (respectively $\mathcal{H}^m_\Z= \{M^m\}/\sim_{\Z H}$).
\end{definition} 

We will show that $\rho_{n}$ is an invariant of $\Q$-homology cobordism.
Since two manifolds that are homology cobordant are necessarily rationally homology cobordant, $\rho_{n}$ is an invariant of homology cobordism.

\begin{theorem}\label{thm:homcob}If $M^{4k-1}_1$ is $\Q$-homology cobordant to
$M^{4k-1}_2$ $(k \geq1)$ then $\rhon(M_1)=\rhon(M_2)$.
\end{theorem}

\begin{proof}Let $W$ be a $4k$-dimensional manifold such that
$\partial W = M_1 \cup \overline{M_2}$, $i_j:M_j \rightarrow W$ 
be the inclusion maps, $E=\pi_1(W)$, and $G_j=\pi_1(M_j)$ for
$j=1,2$.  Since $(i_j)_{\ast}:H_k(M_j;\Q) \rightarrow H_k(W;\Q)$ is 
an isomorphism for $k=1,2$, 
$(i_j)_{\ast}:H_1(G_j;\Q) \rightarrow H_1(E;\Q)$ is an isomorphism
and 
$(i_j)_{\ast}:H_2(G_j;\Q) \rightarrow H_2(E;\Q)$ is surjective.
 Hence by Theorem~4.1 of \cite{CH}, for each $n \geq 0$, the inclusion
maps induce monomorphisms 
\[(i_j)_\ast : \frac{G_j}{(G_j)_H^{(n+1)}} \into
\frac{E}{E_H^{(n+1)}}.
 \]
 Let $\Gamma_n = E/E_H^{(n+1)}$ then we have coefficient systems
$(\beta_j)_n:G_j \rightarrow \Gamma_n$ defined by $(i_j)_\ast \circ
(\phi_j)_n$ where 
 $(\phi_j)_n:G_j \onto G_j/(G_j)_H^{(n+1)}$ is the quotient map.  By
Remark~\ref{fiveremarks}~(\ref{gamma_induction}), 
we have 
 $\rhon(M_j)=\rho(M_j, G_j \rightarrow \Gamma_n)$. 
 Therefore
 \[\rho_{n}(M_1)-\rho_{n}(M_2)=\rho(\partial
W,\Gamma_n)=\st(W,\Gamma_n)-\sigma(W).
 \] 

 To finish the proof, we show that
 $\st(W,\Gamma_n)=\sigma(W)=0$.  
 Since $(i_1)_\ast : H_2(M_1;\Q) \twoheadrightarrow H_2(W;\Q)$ is
surjective,  the second homology of $W$ comes from the boundary.
Thus the 
 intersection of any two classes is $H_2(W;\Q)$ is zero.  In
particular $\sigma(W)=0$.  Let $\mathcal{K}_n$ be the classical right
ring
 of quotients of $\Z\Gamma_n$. Since $H_i(W,M_1;\Q)=0$ for $i=0,1,2$
then by Proposition 2.10 of \cite{COT}, $H_i(W,M_1;\mathcal{K}_n)=0$
for $i=0,1,2$.  By Remark~\ref{fiveremarks}~(\ref{remark:from_bound}), 
$\st(W,\Gamma_n)=0$.
\end{proof}

Hence, for each $n\geq0$ and $m=4k-1$ with $k \geq 1$ we have have a map
\[\rhon:\mathcal{H}^m_\Q \rightarrow \R.\]

\noindent Note that if $L_{1}$ and $L_{2}$ are concordant links then their $0$-surgeries are homology cobordant  hence $\rho_{n}(L_{1})=\rho_{n}(L_{2})$.  
Moreover, $\rho_{n}$ of the trivial link is $0$ as in Example~\ref{example:TrivLink}.  

\begin{corollary}\label{cor:ConcordanceInvt} For each $n\geq 0$, $\rho_{n}$ is a concordance invariant of links and is $0$ for slice links.
\end{corollary}

\noindent We will use $\rho_{n}$ to further investigate the structure of concordance classes of links in Section~\ref{section:filtrations}. 

\section{Non-triviality of $\rho_n$}\label{section:non-triviality}

We will show that the $\rho_{n}$ are highly non-trivial.  To do this we will show 
that the image of $\rhon:\mathcal{H}^3_\Q \rightarrow
\R$ in $\R$ is dense and is an infinitely generated subset 
 of $\R$.  Before we can do this, we must define a family of examples of $3$-manifolds on which we can calculate $\rho_{n}$. 

\subsection{Examples: Genetic Modification}
We describe a procedure wherein one starts with a $3$-manifold $M$ (respectively a link $L$ in $S^{3}$) and ``infects'' $M$ (respectively $L$) along a curve
 $\eta$ in $M$ (respectively $S^{3}-L$) with a knot $K$ in $S^3$ to obtain a new $3$-manifold $M(\eta,K)$ with the same homology as $M$
 (respectively link $L(\eta,K)$ in $S^{3}$).  This construction is a specific type of satellite constructions which has been dubbed \emph{genetic infection} (see Section 3 of \cite{COT1}).

We first describe the construction for a general $3$-manifold.  Let $M$ be a compact, connected, oriented $3$-manifold, $\eta$ be a curve embedded in $M$, and $K$ be a knot in $S^{3}$.  
Denote by $N(\eta)$ and $N(K)$ a tubular neighborhood of $\eta$ in $M$ and $K$ in $S^{3}$ respectively.
Let $\mu_\eta$ and $\mu_K$ be the meridians of $\eta$ and $K$
respectively, and let $l_\eta$ and $l_K$ be the longitudes of $\eta$ and $K$
respectively.  Note that if $\eta$ is not nullhomologous then the longitude of $\eta$ is not well-defined.  In this case, we choose $l_{\eta}$ to be an embedded curve on $N(\eta)$ that is isotopic to $\eta$ in $M$ and intersects $\mu_{\eta}$
geometrically once.
Define 
\begin{equation}\label{examples}M(\eta,K)= (M- \text{N}(\eta)) \cup_f (S^3- \text{N}(K))\end{equation} where
$f: \partial(S^3-\text{N}(K)) \rightarrow \partial(M-\text{N}(K))$ is defined by
$f_\ast(\mu_K)=l_\eta^{-1}$ and $f_\ast(l_K)=\mu_\eta$.   If $\eta$ is not nullhomologous, then there is a choice of longitude for $\eta$ and the homeomorphism type of $M(\eta,K)$ will depend on this choice.  
Since $H_{\ast}(S^{3} - \text{N}(K))$ is independent of $K$, an easy Mayer-Vietoris argument shows that $M$ and $M(\eta,K)$ have isomorphic 
homology groups.
  
Now, consider the case when $M=S^{3}-\text{N}(L)$ where $L$ is an $m$-component link in $S^{3}$ and $\eta$ is a curve in $S^{3}-\text{N}(L) \subset S^{3}$.  Even in the case that 
$\eta$ is not nullhomologous in $H_{1}(S^{3}-N(L))$ there is still a well defined longitude $l_{\eta}$ for $\eta$ since $\eta$ is nullhomologous in $S^{3}$.  
By choosing this longitude, we have a well-defined manifold $M(\eta,K)$.
We now further assume that $\eta$ bounds an embedded disk $D$ in $S^{3}$.
It is
  well known that in this case, $M(\eta,K)$ is homeomorphic to $S^{3}-\text{N}(L(\eta,K))$ where $L(\eta,K)$ is another m-component link in $S^{3}$. 
  Moreover, one can check that $L(\eta,K)$ can be obtained by the following construction.  
Seize the collection of parallel strands of $L$ that pass through the disk $D$ in one hand, just as you might grab some hair in preparation for braiding. Then, treating the collection as a single fat strand, tie it into the knot $K$.
Note that in the special case that  $\eta$ is a meridian of the $i^{th}$  component $L_{i}$ of $L$ then $L(\eta,K)$ is the link obtained adding a local knot $K$ to $L_{i}$.  

\begin{remark} If $L$ is a boundary link in $S^{3}$ then $L(\eta,K)$  is also a boundary link in $S^{3}$.  Hence 
$T(\eta,K)$ is always a boundary link where $T$ is the trivial link.
\end{remark}

\begin{remark} Let $M_{L}$ be the  result of performing $0$-framed surgery on a link $L$ in $S^{3}$ with all linking numbers $0$ and let $\eta$ be a curve in 
$S^{3}-\text{N}(L)\subset M_{L}$ that bounds an embedded disk in $S^{3}$.  Then $M_{L}(\eta,K)=M_{L(\eta,K)}$.  
\end{remark}

\begin{figure}[ht]
\begin{picture}(100,115)
\put(0,0){\includegraphics{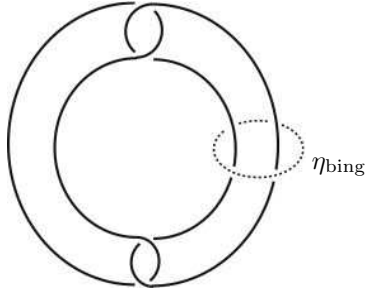}}
\put(115,45){$\eta_{\text{bing}}$}
\end{picture}
\caption{$\eta_{\text{bing}}\in S^{3}-\{\text{trivial link}\}$ }
\label{fig:bingeta} 
\end{figure}

\begin{figure}[ht]
\begin{picture}(100,115)
\put(0,0){\includegraphics{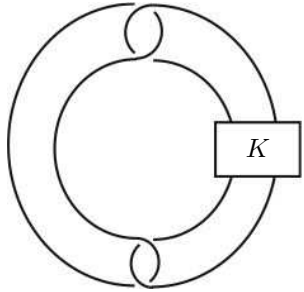}}
\put(90,50){$K$}
\end{picture}
\caption{Bing double of K}
\label{fig:bing}
\end{figure}

\begin{example}
[\textbf{Iterated Bing doubles of $K$}]  Let $T$ be the trivial link with $2$ components and let $\eta_{\text{bing}}$ be the curve in Figure~\ref{fig:bingeta}. 
 Then $\eta_{\text{bing}}$ bounds a disk in $S^{3}$.
$L(\eta,K)$ is the link in Figure~\ref{fig:bing} and is more commonly 
known as the (untwisted) Bing double of $K$, $\text{BD}(K)$.   
We note that $\eta_{\text{bing}} \in F^{(1)} - F^{(2)}$ where $F=\pi_{1}(S^{3}-\text{N}(T))$.
Moreover, any (untwisted) iterated Bing double of $K$ can be obtained as $T(\eta,K)$ where $T$ is a trivial link with $m\geq 2$ components  and $\eta \in F^{(n)}-F^{(n+1)}$ for some $n\geq 1$.
\end{example}

We now construct a cobordism $C=C(\eta,K,W)$ between $M$ and $M(\eta,K)$ for which the inclusion maps 
will behave nicely modulo the torsion-free derived series of their respective fundamental groups.
Recall that $M_K$, 0-surgery on $K$ in $S^3$, is defined as $M_K=
(S^3-\text{N}(K)) \cup_q \text{ST}$ where 
$\text{ST} = S^1 \times D^2$ and $q: \partial(\text{ST}) \rightarrow
\partial(S^3-\text{N}(K))$ is defined so that 
$q_\ast(\{pt\} \times \partial D^2)=l_K$.
By \cite[p. 118]{COT1}, $M_K$ bounds a compact, oriented $4$-manifold $W$ such that 
$\sigma(W)=0$ and 
$\pi_1(W)\cong\Z$, generated by $\iota_\ast(\mu_K)$ where $\iota:M_K
\rightarrow W$ is the inclusion map.
We have an exact sequence $\pi_2(W) \rightarrow H_2(W) \rightarrow
H_2(\pi_1(W))$. Since $\pi_1(W)\cong \Z$, the last term is zero. Hence $H_{2}(W)=\text{image}(\pi_{2}(W) \rightarrow H_{2}(W))$.
Using the $4$-manifold W, we can construct a cobordism $C=C(\eta,K,W)$
between $M$ and $M(\eta,K)$ as follows.  Glue $W$, as above, to
$M\times I$ by identifying $\text{ST}\subset M_K=\partial W$
to $\text{N}(\eta)=D^2 \times S^1\subset M\times \{1\}$ so that
$l_\eta^{-1}$ is identified with $\mu_K$ and $\mu_\eta$ is identified
to $l_K$.
It follows that $\partial C =M \sqcup \overline{M(\eta,K)}$.  Let
$i:M\rightarrow C$ and $j:M(\eta,K) \rightarrow C$ be the inclusion
maps. 
 
 We will use this cobordism to show that the difference between $\rho_{i}(M)$ and $\rho_{i}(M(\eta,K))$ 
 depends only on  $\rho_{0}(K)$ and $\max\{n \hspace{2pt} | \hspace{2pt} \eta \in (\pi_{1}(M))_{\sss H}^{\sss (n)}\}$ (see Theorem~\ref{diff_rho} below).
 We begin with some algebraic lemmas that will be employed in the proof of Theorem~\ref{diff_rho}.
 
 \begin{lemma}\label{istar_iso}$i_\ast:\pi_1(M)\rightarrow \pi_1(C)$
is an isomorphism.
 \end{lemma}
 \begin{proof} $\pi_1(W)\cong \Z$ which is generated by
$\iota_\ast(\mu_K)$.  Moreover, $\mu_K$ generates $\pi_1(W \cap
(M\times I))=\pi_1(\eta\times D^2)$.   Hence, by the Seifert
Van-Kampen Theorem, 
 $\pi_1(M\times I) \rightarrow \pi_1(C)$, induced by the inclusion
map, is an isomorphism. Therefore $i_\ast$ is an isomorphism.
 \end{proof}
 
\noindent  In particular, the inclusion map induces isomorphisms between $\pi_{1}(M)/\pi_{1}(M)_{\sss H}^{\sss {(i)}}$ and 
 $\pi_{1}(C)/\pi_{1}(C)_{\sss H}^{\sss {(i)}}$ for all $i$.  In order to show that $j$ induces an isomorphism between 
 $\pi_{1}(M(\eta,K))/\pi_{1}(M(\eta,K))_{\sss H}^{\sss {(i)}}$ and $\pi_{1}(C)/\pi_{1}(C)_{\sss H}^{\sss {(i)}}$ for all $i$, we
appeal to Theorem~\ref{CH_thm} and Proposition~\ref{CHprop}. The next two lemmas will 
 guarantee that the hypotheses of Proposition~\ref{CHprop} are satisfied.   
 
\begin{lemma} \label{GE_epi}$j_\ast:\pi_1(M(\eta,K)) \rightarrow
\pi_1(C)$ is an epimorphism. 
\end{lemma}
\begin{proof} Let $\alpha$ be a curve in $\pi_1(C)$.  By
Lemma~\ref{istar_iso}, $i_\ast$ is an isomorphism. Hence $\alpha$ can
be represented
by a curve in $M$.  Moreover, by general position, we can assume
$\alpha$ misses $\text{N}(\eta)$.  Push $\alpha$ into $M\times \{1\}$ to get
a 
curve $\beta$ in $(M\times \{1\}) - (\text{N}(\eta)\times\{1\}) \subset
M(\eta,K)$ such that $j_\ast(\beta)=\alpha$. 
\end{proof}

\begin{lemma}\label{lemma3}
 $j_\ast :\pi_{1}(M(\eta,K)) \rightarrow \pi_{1}(C)$ induces an isomorphism on $H_1(-;\Z)$ and
an epimorphism on $H_2(-;\Z)$.
\end{lemma}
\begin{proof}
Let $G=\pi_1(M(\eta,K))$ and $E=\pi_1(C)$.
Since $j$ induces an epimorphism on $\pi_1$, $j$ induces an
epimorphism on $H_1(-)$.  
The inclusion $l:
M-\text{N}(\eta) \rightarrow M(\eta,K)$ induces a epimorphism on $H_{1}(-;\Z)$.
Moreover, the inclusion $l^\prime : M-\text{N}(\eta) \rightarrow M$ induces
an epimorphism on $H_{1}(-;\Z)$ and
the kernel of $l^\prime_\ast$ is the 
subgroup generated by $\mu_\eta \in H_1(M-\text{N}(\eta))$.  

Consider the following  commutative diagram where all of the maps are
induced
from inclusion maps: 
\[
\begin{diagram}
H_1(M-\text{N}(\eta)) & \rOnto^{l^\prime_\ast} & H_1(M) \\
\dOnto_{l_\ast} & & \dTo_{i_\ast}^\cong\\
H_1(M(\eta,K)) & \rOnto^{j_\ast} & H_1(C).
\end{diagram}
\]
Suppose $\alpha \in H_1(M(\eta,K))$ and $j_\ast(\alpha)=0$.  Since
$l_\ast$ is surjective, there exists $\gamma \in H_1(M-\text{N}(\eta))$ such
that
$l_\ast(\gamma)=\alpha$.  Hence
$i_\ast(l^\prime_\ast(\gamma))=j_\ast(l_\ast(\gamma))=0$.  Since
$i_\ast$ is injective, $\gamma \in \ker(l^\prime_\ast)$ so
$\gamma$ is a multiple of $\mu_\eta$.   Howeover,
$l_\ast(\mu_\eta)=0$ since in $M(\eta,K)$, $\mu_\eta$ is identified
with $l_K$, which bounds 
a surface in $S^3-K \subset M(\eta,K)$.  Therefore $j_\ast:
H_1(M(\eta,K))\rightarrow H_1(C)$ is a monomorphism.  Since
$H_1(M(\eta,K))=H_1(G)$ (similarily for $E$),
$j_\ast: H_1(G)\rightarrow H_1(E)$ is an isomorphism.

Let $\alpha\in H_2(E)$, then there exists $\beta\in H_2(C)$ such that
$k(\beta)=\alpha$ where $k:H_2(C) \rightarrow H_2(E)$
is the map in the exact sequence
\[\pi_2(C) \rightarrow H_2(C) \xrightarrow{k}  H_2(E) \rightarrow 0.
\]
  By a Mayer-Vietoris argument, we see that 
$H_2(C)\cong H_2(W) \oplus H_2(M\times I)$ (in the obvious way).
Since $H_2(W)= \text{im}(\pi_2(W)\rightarrow H_{2}(W))$
and any element of $\pi_2(W)$ goes
to zero under $k$, we can assume that $\beta \in H_2(M\times I)\oplus
\{0\} \cong H_2(M)$.  Let $S_\beta$ be a surface in $M$ representing
$\beta$.  We will construct a surface $S^\prime_\beta \subset
M(\eta,K)$ such that
$k(j_\ast([S^\prime_\beta]))=k(i_{\ast}([S_\beta]))=\alpha$.  This will
complete the proof
since $k \circ j_\ast = j_\ast \circ k^\prime$ where $k^\prime:
H_2(M(\eta,K)) \onto H_2(G)$.   

We can assume that $S_\beta$ intesects 
$\text{N}(\eta)$ in finitely many disks $D_i$.   To construct
$S_\beta^\prime$, remove each of these disks and replace them with a
copy of a chosen seifert surface for $K$ in $S^3-K \subset M(\eta,K)$
oriented according the signed intersection
of $\eta$ with $S_\beta$. Let $F$ be the surface in $W \subset C$
obtained by gluing a disk in $(M\times \{1\}) \cap W$ whose boundary
is $\mu_\eta$ with a copy of the chosen seifert surface for $K$.
Since $H_2(\pi_1(W))=0$, 
$k([F])$ is trivial in $H_2(E)$.   Moreover
$[S^\prime_\beta]=m[F]+[S_\beta]$ in $H_{2}(C)$ for some $m$ so
$k(j_\ast([S^\prime_\beta]))=k(i_{\ast}([S_\beta]))$.
\end{proof}

Note that $M_{K}=\partial W \subset C$.  For each $i\geq 0$, let $\tau_{i}:\pi_1(M_K) \rightarrow
\pi_{1}(C)/{\pi_{1}(C)}^{\sss (i+1)}_{\sss H}$ be the composition of the map induced by inclusion $\pi_{1}(M_{K}) \rightarrow \pi_{1}(C)$ and 
the quotient map $\pi_{1}(C) \rightarrow \pi_{1}(C)/{\pi_{1}(C)}^{\sss (i+1)}_{\sss H}$.

 \begin{lemma}\label{imagetau}If $\tau_i: \pi_1(M_K) \rightarrow \pi_{1}(C)/{\pi_{1}(C)}^{\sss (i+1)}_{\sss H}$ is the homomorphism as described above then
$$Im( \tau_i)\cong 
\left\{
\begin{array}{ll}
    \{1\} & 0 \leq i\leq  n-1;\\
    \Z & i \geq n.\\
\end{array}
\right.    $$   
\end{lemma}
\begin{proof}  Let $E=\pi_{1}(C)$.
Recall that  $\pi_1(M_K) \cong \pi_1(S^3-K)/\left<l_K\right>$ where
$l_K$ is the longitude of $K$.  Hence every element of $\pi_1(M_K)$
can be
represented by $\alpha = \prod_i g_i \mu_K g_i^{-1}$
where $g_i \in \pi_1(S^3-K)$ and $\mu_K$ is a fixed meridian 
of $K$.  Let $\tau : \pi_1(M_K) \rightarrow E$ be induced by the
inclusion of $M_K$ into $C$.  
Since $\mu_K$ is identified to $l_\eta^{-1}$ in $E$, we see that
$\tau(\alpha)=\prod_i \tau(g_i) i_\ast(l_\eta)^{-1} \tau(g_i)^{-1}$.  
Moreover, since $l_\eta \in P_{^H}^{_{(n)}}$ and $P/P_{^H}^{_{(i+1)}}
\cong E/E_{^H}^{_{(i+1)}}$ for all $i$, $\tau(\alpha) \in
E_{^H}^{_{(i+1)}}$ for $0 \leq i \leq n-1$. 
Therefore the image of $\tau_i$ is trivial for $0 \leq i \leq n-1$.

We will prove that the image of $\tau_i$ is $\Z$ for $i\geq n$ by
induction on $i$.  First we prove this is true for $i=n$.
Since $\tau(\pi_1(M_K)) \subset E_{^H}^{_{(n)}}$,
$\tau([\pi_1(M_K),\pi_1(M_K)])\subset 
[E_{^H}^{_{(n)}},E_{^H}^{_{(n)}}] \subset E_{^H}^{_{(n+1)}}$.  Hence
we have a well-defined map 
$$\overline{\tau}: \Z \cong \pi_1(M_K)/\pi_1(M_K)^{(1)} \rightarrow
E/E_{^H}^{_{(n+1)}}.$$  Since $l_\eta \not\in P_{^H}^{_{(n+1)}}$ and
$P/P_{^H}^{_{(n+1)}} \cong E/E_{^H}^{_{(n+1)}}$, 
$\tau(\mu_K)= i_\ast(l_\eta^{-1}) \not\in E_{^H}^{_{(n+1)}}$.
Therefore $\overline{\tau}$ is non-trivial.  
Moreover, since $\tau(\pi_1(M_K)) \subset E_{^H}^{_{(n)}}$, the image
of $\overline{\tau}$ is contained is
$E_{^H}^{_{(n)}}/E_{^H}^{_{(n+1)}}$ which is $\Z$-torsion free.  
It follows that the image of $\tau_n$ is isomorphic to $\Z$.

To finish the induction, assume that the image of $\tau_i$ is
isomorphic to $\Z$ for some $i\geq n$.  
Let $A=\pi_1(M_K)$; then 
by Example~\ref{ex:knot}, 
$A/A_{^H}^{_{(j+1)}} \cong
\Z$ for all $j \geq 0$. 
By assumption, $\tau$ induces a monomorphism $\tau_\ast :
A/A_{^H}^{_{(i+1)}} \hookrightarrow E/E_{^H}^{_{(i+1)}}.$     
Thus, 
by Proposition~\ref{functoriality}, 
$\tau$ induces a map $\tau_\ast:
A/A_{^H}^{_{(i+2)}} \rightarrow E/E_{^H}^{_{(i+2)}}$.  
Since  $A_{^H}^{_{(i+1)}}=A_{^H}^{_{(i+2)}}$, the map $\tau_\ast:
A/A_{^H}^{_{(i+2)}} \rightarrow E/E_{^H}^{_{(i+2)}}$ is a
monomorphism.  
Therefore the image of $\tau_{i+1}$ is $\Z$.    
\end{proof}

\begin{theorem}\label{diff_rho}Let $M(\eta,K)$ be as defined in $(\ref{examples})$
and $P=\pi_1(M)$.  If $\eta \in P_H^{(n)} - P_H^{(n+1)}$ for some
$n\geq 0$ then 
$$\rho_i(M(\eta,K))-\rho_i(M)= 
\left\{
\begin{array}{ll}
    0 & 0 \leq i\leq  n-1;\\
    \rho_0(K) & i \geq n.\\
\end{array}
\right.    
$$   
\end{theorem}

\noindent Before proving Theorem~\ref{diff_rho}, we establish two easy corollaries.

\begin{corollary} Let $\eta$ be a embedded curve in $\#_{i=1}^{m} S^{1}\times S^{2}$ that is non-trivial in $F=\pi_{1}(\#_{i=1}^{m} S^{1}\times S^{2})$.  
If $K$ is a knot with $\rho_{0}(K)\neq 0$ then $M(\eta,K)$ is not $\Q$-homology cobordant to $\#_{i=1}^{m} S^{1}\times S^{2}$.
\end{corollary}
\begin{proof} Since $F$ is a free group, $F^{\sss (k)}_{\sss H} = F^{\sss (k)}$ for all $k\geq 0$ and $F^{\sss (\omega)}=\{1\}$. Therefore $\eta \in F^{\sss (n)}-F^{\sss (n+1)}$ for some $n\geq 0$.  
By Theorem~\ref{diff_rho}, $\rho_{n}(M(\eta,K)) =\rho_{0}(K) + \rho_{n}(\#_{i=1}^{m} S^{1}\times S^{2})=\rho_{0}(K)$. Thus, if $K$ is a knot with $\rho_{0}(K)\neq 0$ then
$M(\eta,K)$ is not $\Q$-homology cobordant to $\#_{i=1}^{m} S^{1}\times S^{2}$.
\end{proof}

\begin{corollary}\label{cor:itbingdouble}Suppose $K$ is a knot in $S^{3}$ with $\rho_{0}(K)\neq 0$.  Then no iterated Bing double of $K$ is slice.
\end{corollary}
\begin{proof}Any iterated Bing double of $K$ can be obtained as 
$T(\eta,K)$ where $\eta$ is a non-trivial commutator of length $m+1$ where $m$ is the number of components of $T(\eta,K)$ \cite{C1}.  Therefore, $\eta$ is a nontrivial
element in $F^{(n)}/F^{(n+1)}$ for some $n\geq 1$ where $F=\pi_{1}(T)$.  By Theorem~\ref{diff_rho}, $\rho_{n}(T(\eta,K))=\rho_{0}(K)\neq 0$.  Hence $T(\eta,K)$
is not slice. 
\end{proof}
\noindent Subsequent to our work, it has been shown that if the Bing double of $K$ is slice then $K$ is algebraically slice \cite{CLR}.

 We will now prove Theorem~\ref{diff_rho}.

\begin{proof}[Proof of Theorem~\ref{diff_rho}]  
Let $W$ be as defined before, $V$ be an open neighborhood of
$M_K=\partial W$ in $W$, $W_V = W-V$, and $C_V=\overline{C-{W_V}}$
where $C=C(\eta,K,W)$ is as described in the discussion preceding Lemma~\ref{istar_iso}.
Then $V$ is 
homeomorphic to $M_K \times [0,1]$ and $C_V$ can be obtained by
gluing $M \times [0,1]$ to $M_K \times [0,1]$ where $\eta \times D^2$ 
is identified to $\text{ST} = D^2 \times S^1$ in $M_K$.  It follows
that $C_V$ has $3$ boundary components; in particular, $\partial C_V =
M \sqcup \overline{M(\eta,K)} \sqcup M_K$.

\begin{figure}[htbp]
\begin{center}
\begin{picture}(165,100)
%\put(64,41){$\theta$}
\put(0,0){\includegraphics{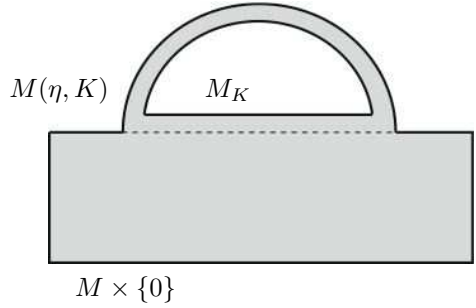}}
\put(10,-12){$M \times \{0\}$}
\put(-15,63){$M(\eta,K)$}
\put(59,63){$M_K$}
\end{picture}
\end{center}
\caption{The $4k$-manifold $C_V$ with $\partial C_V = M \sqcup
\overline{M(\eta,K)} \sqcup M_K$}\label{fig:3boundcomps}
\end{figure}

Let $G=\pi_1(M(\eta,K))$, $E=\pi_1(C)$, and 
$\Gamma_i=E/E_{^H}^{_{(i+1)}}$ for $i\geq0$.  Since $C_V \subset
C$, there is an obvious coefficient system $\pi_1(C_V) \rightarrow
\Gamma_i$ (similarly for $M$, $M(\eta,K)$, and $M_K$).  By
Lemmas~\ref{istar_iso}, 
\ref{GE_epi} and \ref{lemma3}, Theorem~\ref{CH_thm}, and Proposition~\ref{CHprop}, $i_\ast
:P/P_{^H}^{_{(i+1)}} \rightarrow \Gamma_i$ and $j_\ast :
G/G_{^H}^{_{(i+1)}} \rightarrow \Gamma_i$ are isomorphisms for all
$i\geq 0$.
Therefore, $\rho_i(M)=\rho(M,\Gamma_i)$ and
$\rho_i(M(\eta,K))=\rho(M(\eta,K),\Gamma_i)$.  Hence 
we have
\begin{equation}\label{rho_sig}\rho_i(M)-\rho_i(M(\eta,K)) +
\rho(M_K,\Gamma_i) = \st(C_V,\Gamma_i)-\sigma(C_V)
\end{equation}
for all $i\geq0$.  

Since the abelianization of $\pi_1(M_K)$ is $\Z$, there is a unique
surjective homomorphism $\pi_1(M_K) \onto \Z$ up to isomorphism.
Therefore, by Lemma~\ref{imagetau} and Lemma~\ref{gamma_induction_for_rho} we have $\rho(M_K,\Gamma_i)=\rho_0(M_K)$
for all $i \geq n$, and $\rho(M_K,\Gamma_i)=0$ for $0\leq i\leq n-1$. 

To finish the proof, we will show that $\st(C_V,\Gamma_i)-\sigma(C_V)=0$
for all $i\geq 0$.  Since $\mathcal{U}\Gamma_i$ is a flat
$\Z\Gamma_i$-module 
it suffices to
show that the map
induced by the inclusion $H_2(\partial C_V;\Z\Gamma_i) \rightarrow
H_2(C_V;\Z\Gamma_i)$ is surjective for $-1 \leq i$ (recall that
$\Gamma_{-1}=\{1\})$.   

Recall that  $C_V=M_K \times I \cup_{\eta \times D^2} M \times I$.  
Consider the Mayer-Vietoris sequence
\begin{align*} \rightarrow H_2(M_K \times I;\Z\Gamma_i) \oplus H_2(M
\times I;\Z\Gamma_i) \rightarrow H_2(C_V;\Z\Gamma_i)  \rightarrow
H_1(\eta \times D^2;\Z\Gamma_i)  \\
 \rightarrow H_1(M\times I;\Z\Gamma_i) \oplus H_1(M_K\times
I;\Z\Gamma_i). 
\end{align*} 
By Lemma~\ref{imagetau}, if $i\leq n-1$ then $M_K$ lifts to the
$\Gamma_i$-cover hence $H_1(M_K \times I;\Z\Gamma_i)\cong \Z\Gamma_i$ is
generated by $\mu_K$. 
Moreover, in this case, $\eta$ lifts to the $\Gamma_i$-cover and
$H_1(\eta \times D^2;\Z\Gamma_i) \rightarrow H_1(M_K;\Z\Gamma_i)$ is
an isomorphism.  Therefore 
$H_2(M_K\times I;\Z\Gamma_i)\oplus H_2(M \times I;\Z\Gamma_i)
\rightarrow H_2(C_V;\Z\Gamma_i)$ is surjective for $i\leq n-1$.  
If $i\geq n$, then the image of $\pi_1(\eta \times D^2)$ in
$\Gamma_i$ is $\Z$ by Lemma~\ref{imagetau}. Hence $H_1(\eta \times
D^2;\Z\Gamma_i)=0$. Thus
 for $i\geq n$ we have that $H_2(M_K\times I;\Z\Gamma_i)\oplus H_2(M
\times I;\Z\Gamma_i) \rightarrow H_2(C_V;\Z\Gamma_i)$ is surjective.
It follows that $H_2(\partial C_V;\Z\Gamma_i) \rightarrow
H_2(C_V;\Z\Gamma_i)$ is surjective for all $i\geq -1$.
\end{proof}
 
\subsection{Non-triviality of Examples}
 
For each knot $K$ in $S^3$, there exists a degree one map
$f_K: S^3-K \rightarrow S^4-T$ (where $T$ is the unknot) that
induces an isomorphism on homology with $\Z$ coefficients and
fixes the boundary.  Hence, there is a degree one map
$\overline{f}_K: M(\eta,K) \rightarrow M(\eta,T)=M$ that induces an
isomorphism on $H(-;\Z)$.   Recall that $\mathcal{H}_{\Q}^m$ is the set of $\Q$-homology cobordism
classes of closed, oriented $m$-dimensional manifolds.  
For a fixed closed, oriented $m$-dimensional manifold $M$,
 we define $\mathcal{H}^m_\Q(M) \subset \mathcal{H}^m_\Q$ as follows:
$[N^\prime] \in \mathcal{H}^m_\Q(M)$ if there exists an $N$
  such that $[N]=[N^\prime] \in \mathcal{H}^m_\Q$ and there exists
a degree one map $f: N \rightarrow M$
 that induces an isomorphism on $H_\ast(-;\Q)$. 
 By Proposition~\ref{prop:connsum}, $\rho_{n}$ is additive under the connected sum 
 of manifolds.  Therefore the images of $\rho_{n}:\mathcal{H}^m_\Q(M) \rightarrow \R$
 and $\rho_{n}:\mathcal{H}^m_\Q \rightarrow \R$ are subgroups of $\R$.

\begin{theorem}\label{gen_dense}Let $M$ be a closed, oriented
$3$-manifold, $G=\pi_1(M)$ and $n \geq 0$. If $\gn/\gnp \neq \{1\}$
then the image
of $\rho_n : \mathcal{H}^3_\Q(M) \rightarrow \R$ is (1) dense in $\R$
and (2) an infinitely generated subgroup of $\R$.
In particular, the image of 
of $\rho_n : \mathcal{H}^3_\Q\rightarrow \R$ is (1) dense in $\R$
and (2) an infinitely generated subgroup of $\R$ .
\end{theorem}
\begin{proof}  Let $-2 < r < 2$.  By \cite[Section~2]{CL}, for all
$\epsilon > 0$, there exists a knot $K_{r,\epsilon}$ in $S^3$ such
that 
$$ \left|\int_{S^1} \sigma_\omega(K_{r,\epsilon}) d\omega - r\right|
< \epsilon.$$ 
Here, $\sigma_\omega(K)$ is the Levine-Tristram signature of $K$ at
$\omega \in S^1$ and the circle is normalized to have length $1$.
By \cite[Lemma 5.4]{COT}, $\rho_0(M_K) = \int_{S^1} \sigma_\omega(K)
d\omega$  for any knot $K$ in $S^3$.  
Let $\eta$ be a curve in $M$ representing an element in $\gn - \gnp$.
Then by Theorem~\ref{diff_rho}, $\rho_n(M(\eta,K))=\rho_0(M_K)$.  
Moreover, by the above remarks, there is a degree one map from
$M(\eta,K)$ to $M$ hence $r \in
\overline{\rho_n(\mathcal{H}^3_\Q(M))}$.

For arbitrary $r \in \mathbb{R}$ there exists a positive integer $m$
such that $r/m \in (-2,2)$.  Let $\epsilon > 0$ and set
$\epsilon^\prime = \epsilon/m$.  
As above, there exists a knot $K_{r/m,\epsilon^\prime}$ such that
$|\rho_0(M_{K_{r/m,\epsilon^\prime}})-r/m| < \epsilon^\prime$.  We
remark that $K_1 \# K_2$ 
can be obtained as $M_{K_1}(\eta,K_2)$ where $\eta$ is a meridian of
$K_1$ in $M_{K_1}$.  Hence by Theorem~\ref{diff_rho}, 
$\rho_0(M_{K_1 \# K_2})=\rho_0(M_{K_1}) + \rho_0(M_{K_2})$
 (this can also be computed directly).  Hence
$\rho_0(M_{mK})=m\rho_0(M_K)$ where $mK$ is the connected sum of $K$
with itself $m$ times.  Therefore
 $$ \left|\rho_0(M_{mK_{r/m,\epsilon^\prime}}) - r \right| = m \left|
\rho_0(M_{K_{r/m,\epsilon^\prime}}) - r/m \right| < m \epsilon^\prime
= \epsilon.$$
As before, it follows that $r \in
\overline{\rho_n(\mathcal{H}^3_\Q(M))}$.  This completes the verification that the image of $\rho_{n}$ is dense in $\R$.

By Proposition~2.6 of \cite{COT1}, there exists an infinite set
$\{J_i | i \in \Z_+\}$ of Arf invariant zero knots such that
$\{\rho_0(M_{J_i})\}$ is linearly
independent over the integers.   Therefore $\{\rho_n(M(\eta,J_i))\}$
is an infinitely generated subgroup of $\R$.  Since $\rho_{n}(\mathcal{H}^{3}_{\Q}(M))$ is an abelian group that contains this subgroup, it is itself infinitely generated.  
\end{proof}

Let $T$ be a trivial link.  
If $L$ is a boundary link then, just as for a knot, there is a
degree one map $S^3-L$ to $S^3-T$ that
fixes the boundary and induces an 
isomorphism on homology.
In particular, there is a degree one map from the $0$-surgery on a boundary
link with $m$ components to the connected sum of $m$ copies of $S^1
\times S^2$ that induces
an isomorphism on homology.  Thus, for each $m\geq 1$, we can
consider the subset $\mathcal{H}^{3,b}_\Q(m) \subset
\mathcal{H}^3_\Q(\#_{i=1}^m S^1 \times S^2)$ defined by 
$[N^\prime] \in \mathcal{H}^{3,b}_\Q$ if $[N^\prime]= [N] \in
\mathcal{H}^3_\Q$ where $N$ is 0-surgery on a boundary link in $S^3$.
As before, for each $n$ and $m$, the image of $\rho_{n}:\mathcal{H}^{3,b}_\Q(m) \rightarrow \R$
is a subgroup of $\R$.

\begin{corollary} \label{dense} For each $n\geq 0$ and $m\geq 2$,
the image of $\rho_n : \mathcal{H}^{3,b}_\Q(m) \rightarrow \mathbb{R}$ is (1) dense in $\mathbb{R}$
and (2) an infinitely generated subgroup of $\R$.   
\end{corollary}
\begin{proof}
Let $M$ be the manifold obtained by performing 0-surgery on the
trivial link with $m\geq 2$ components.  Then $M=\#_{i=1}^m S^1
\times S^2$ and 
$\pi_1(M)\cong F(m)$. 
Since 
$F^{_{(n)}}_{^H}/F^{_{(n+1)}}_{^H}=F^{(n)}/F^{(n+1)}$
is 
non-trivial for $n\geq 0$, there exists an $\eta \in F^{_{(n)}}_{^H}
-F^{_{(n+1)}}_{^H}$.  We can assume that $\eta \in \pi_1(S^3-T)$
where $T$ is the trivial link in $S^3$ since 
the inclusion $S^3-T \rightarrow M$ induces an isomorphism on
$\pi_1$. Moreover, we can alter $\eta$ by a homotopy  in $S^3-T$ to
obtain a curve $\eta^\prime$ that bounds a disk in $S^3$.  
Thus, as we showed in the proof of Theorem~\ref{gen_dense}, the image of $\rho_n : \{M(\eta^\prime,K)| K \text{ is a knot
in }S^3\} \subset \mathcal{H}^{3,b}_\Q(m) \rightarrow \R$ is dense
in $\R$ and is an infinitely generated subgroup of $\R$.
\end{proof}

We will now show that the $\rho_n$ are independent functions.  For each
$m\geq 2$, let $V_m = \{f: \mathcal{H}^{3,b}_\Q(m) \rightarrow
\R\}$, be the vector space of functions from the set
$\mathcal{H}^{3,b}_\Q(m)$ to $\R$ and $V=\{f:\mathcal{H}^3_\Q
\rightarrow \R\}$ 
be the vector space of functions from $\mathcal{H}^3_\Q$ to $\R$.

\begin{theorem} \label{thm:indep} For each $m \geq 2$, $\{\rho_n \}_{n=0}^{\infty}$
is a linearly independent subset of $V_m$.  In particular,
$\{\rho_n\}_{n=0}^{\infty}$ is a linearly 
independent subset of $V$.
\end{theorem}
\begin{proof} Let $\alpha = r_1 \rho_{i_1} + \cdots + r_k \rho_{i_k}$
where $r_i$ are non-zero real numbers, $\rho_{i_j} : \mathcal{H}^{3,b}_\Q(m)
\rightarrow \R$  and $i_{j} < i_{j^\prime}$ if $j < j^\prime$. 
Suppose $\alpha=0$.  Let $M$ be 0-surgery on the $m$-component
trivial link $T$ in $S^3$ and $F=\pi_1(M)$.  As in the proof of 
Corollary~\ref{dense} above, for each $n\geq 0$, there is a curve $\eta_n$
in $S^3-T$ such that $\eta_n$ bounds a disk in $S^3$ and $\eta_n \in
F^{_{(n)}}_{^H} -F^{_{(n+1)}}_{^H}$.
For each $n\geq 0$, let $M_n=M(\eta_n,K)$ for some $K$ with
$\rho_0(M_K)\neq 0$ (for example, let $K$ be the right handed
trefoil).  
By the remarks above Corollary~\ref{dense}, $M_n \subset
\mathcal{H}^{3,b}_\Q(m)$.  By Theorem~\ref{diff_rho},
$\rho_i(M_n)=0$ for $i\leq n-1$ and $\rho_n(M_n)\neq 0$. 
Hence $0=\alpha(M_{i_k})=r_k \rho_{i_k}(M_{i_k})$.  Since
$\rho_{i_k}(M_{i_k})\neq 0$, $r_k=0$.  This is a contradiction.
\end{proof}

\section{The Grope and ($n$)-Solvable Filtrations}\label{section:filtrations}

We now investigate the grope and $(n)$-solvable
filtrations of the string link concordance group (with $m \geq 2$
strands) first defined for knots by
 T. Cochran, K. Orr and P. Teichner in \cite{COT}.  We will show that the function 
 $\rho_{n}$ is a homomorphism on the subgroup of boundary links and vanishes for $(n+1)$-solvable links.  
Using this, we will show that each of the successive quotients of the $(n)$-solvable filtration of the string link concordance
group contains an infinitely generated subgroup (even modulo local knotting).  We will also show that a similar statement 
holds for the grope filtration of the string link concordance group.   The reason that we study string links instead of ordinary links is that 
the connected sum operation for ordinary links is not well-defined.  Thus, in order to have a group structure on the set of links up to concordance, we must 
use string links.

We begin by recalling some definitions.  Recall that an \textbf{$m$-component string link} (sometimes
called an m-component disk link \cite{L,LD})
is a  locally flat embedding $f: \sqcup_{m} I \rightarrow D^{3} $ of $m$ oriented, ordered copies of $I$ in $D^{3}$ that is 
transverse to the boundary and such that $f|{\sqcup_{m} \partial I}$ is the standard m-component trivial $0$-link in $S^{2}$, $j_{0}: \sqcup_{m} \partial I \rightarrow S^{2}$. 
Two $m$-component string links $f,g$ are \textbf{concordant} if there exists a locally flat embedding $F:\sqcup_{m} I\times I \rightarrow D^{3}\times I $
that is transverse to the boundary and such that $F|{\sqcup_{m} I\times \{0\}}=f, F|{\sqcup_{m} I\times \{1\}}=g$ and $F|\sqcup_{m} \partial I \times I= j_{0} \times \text{id}_{I} $.
The concordance classes of $m$-component string links forms a group under stacking (see Figure~\ref{fig:stack}) which we denote
by $C(m)$ (see \cite{LD} for more details).  In the literature this group is often denoted
$C(m,1)$ or $CSL(m)$.   This group is known to be non-abelian when $m\geq 2$ 
\cite{LD} and abelian when
$m=1$.  Let $\SB(m)$ be the subgroup of \emph{boundary} disk links in $C(m)$.

\begin{figure}
\begin{picture}(167,138)
\put(0,0){\includegraphics{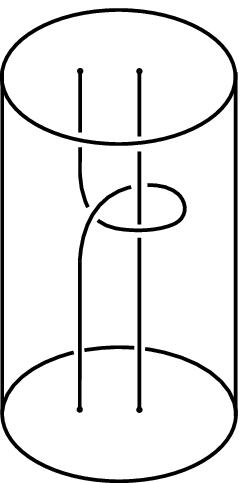}}
\put(-25,70){$L=$}
\put(109,20){\includegraphics{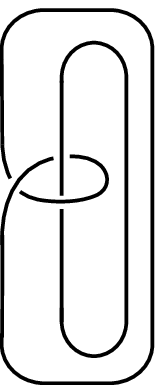}}
\put(160,70){$=\hat{L}$}
\end{picture}
\caption{The string link $L$ and its closure $\hat{L}$}\label{fig:stringlink}
\end{figure}

 If $L$ is a string link then the closure of $L$, denoted by
$\wh{L}$,  is the oriented, ordered $m$-component link in $S^{3}$ obtained by adjoining to its boundary the standard trivial string link.
 A string link is equipped with a well-defined set of meridians that we denote by $\mu_{1},\dots,\mu_{m}$.

\subsection{$(n)$-solvable filtration of $C(m)$}In \cite[Definition 8.5, p. 500]{COT} and \cite[Definition 8.7, p. 503]{COT},  Cochran, Orr, and Teichner define
 an $(n)$-solvable knot and link where  $n \in {\frac{1}{2}}\mathbb{N}_{0}$. 
We recall the definitions here.
We first recall the definition of an $(n)$-Lagrangian.   Let $M$ be a fixed closed, oriented $3$-manifold.  An \textbf{$H_{1}$-bordism}
 is a $4$-dimensional spin manifold $W$ with boundary $M$ such that the inclusion map induces an isomorphism $H_{1}(M)\cong H_{1}(W)$.
 For any $4$-manifold $W$, let $W^{(n)}$ denote the regular covering of $W$ that corresponds to the $n^{th}$ term of the derived series of $\pi_{1}(W)$.  
 For each $n\geq 0$, one can define the quadratic forms $\lambda_{n}, \mu_{n}$ on $H_{2}(W^{(n)})$ in terms of equivariant intersection and self-intersection
 numbers of surfaces in $W$ that lift to the cover $W^{(n)}$.  If $F$ is a closed, oriented, immersed, and based surface in $W$ that lifts to $W^{(n)}$ then $F$ is called an \textbf{(n)-surface}.
See \cite[Chapter 7, pp. 493--496]{COT} for the 
 precise definitions.  
 
 \begin{definition}\label{nlagrangian} Let $W$ be an $H_{1}$-bordism such that $\lambda_{0}$ is a hyperbolic form.  
 \begin{enumerate}
 \item A \textbf{Lagrangian} for $\lambda_{0}$ is a direct summand of $H_{2}(W)$ of half rank on which 
 $\lambda_{0}$ vanishes.
 \item An \textbf{(n)-Lagrangian} is a submodule $L\subset H_{2}(W^{(n)})$ on which $\lambda_{n}$ and $\mu_{n}$ vanish
 and which maps onto a Lagrangian of the hyperbolic form $\lambda_{0}$ on $H_{2}(W)$.
 \item Let $k\leq n$.  We say that an (n)-Lagrangian $L$ \textbf{admits (k)-duals} if $L$ is generated by (n)-surfaces $l_{1},\dots, l_{g}$
 and there are (k)-surfaces $d_{1},\dots, d_{g}$ such that $H_{2}(W)$ has rank $2g$ and 
 \[ \lambda_{k}(l_{i},d_{j})=\delta_{ij}.\]
 \end{enumerate}
 \end{definition}

We now define $(n)$-solvability for a $3$-manifold or link in $S^{3}$.

\begin{definition}  Let $M$ be a closed, oriented $3$-manifold and $n\in \mathbb{N}$.  $M$ is called
\textbf{(0)-solvable} if it bounds an $H_{1}$-bordism $W$ such that $(H_{2}(W),\lambda_{0})$ is
hyperbolic.
 $M$ is called \textbf{(n)-solvable} for $n>0$ if there is an $H_{1}$-bordism $W$ for $M$that
contains an
(n)-Lagrangian with (n)-duals. The $4$-manifold $W$ is called an \textbf{$(n)$-solution} for $M$.  A link $L$ in $S^{3}$ is said to be \textbf{$(n)$-solvable} if the zero surgery on $L$ is $(n)$-solvable.
A string link $L \in C(m)$ is said to be \textbf{$(n)$-solvable}
if $\widehat{L}$ is $(n)$-solvable.
\end{definition}

\begin{definition}Let $M$ be a closed, oriented $3$-manifold and $n\in\mathbb{N}_{0}$.  
M is said to be \textbf{(n.5)-solvable} if there is an $H_{1}$-bordism for $M$ that contains an $(n+1)$-Lagrangian with 
$(n)$-duals.  The $4$-manifold is called an \textbf{$(n.5)$-solution}.  A link $L$ in $S^{3}$ is said to be \textbf{$(n.5)$-solvable} if the zero surgery on $L$ is $(n.5)$-solvable.  A string link $L \in C(m)$ is said to be \textbf{$(n.5)$-solvable}
if $\widehat{L}$ is $(n.5)$-solvable.
\end{definition}
 
 Hence,  for each $m\geq 1$, we can define a filtration of the string link concordance group
 \[ \cdots \subset \mathcal{F}_{(n.5)}^{m} \subset \mathcal{F}_{(n)}^{m} \subset \cdots \subset
\mathcal{F}_{(0.5)}^{m} \subset \mathcal{F}_{(0)}^{m} \subset C(m)\]
by setting $\mathcal{F}^{m}_{(n)}$ to be the set of (n)-solvable $L\in C(m)$ for $n\in
\frac{1}{2}\mathbb{N}_{0}$.  It is easy to check that $\mathcal{F}_{(n)}^{m}$ is a subgroup of
$C(m)$ and in fact is a normal subgroup of $C(m)$ for each $m\in \mathbb{N}$ and $n\in
\frac{1}{2}\mathbb{N}_{0}$.

When $m=1$, $C(1)$ is the concordance group of knots which is an abelian group.  It was shown in \cite{COT}
that $C/\mathcal{F}^{1}_{(0)}\cong \Z_{2}$ given by the Arf invariant and $C/\mathcal{F}^{1}_{(0.5)}$
is J. P. Levine's algebraic concordance group which is isomorphic to  $\Z^{\infty}\oplus
\Z^{\infty}_{2}\oplus \Z_{4}^{\infty}$ \cite{Le}.  It was also shown in \cite{COT,COT1} that
$\mathcal{F}^{1}_{(n)}
/\mathcal{F}^{1}_{(n.5)}$ for $n=1,2$ has infinite rank.  Moreover, in \cite{CT}, 
Cochran and Teichner showed that 
that for each $n\in \mathbb{N}$, $\mathcal{F}^{1}_{(n)} /\mathcal{F}^{1}_{(n.5)}$ is of
infinite order.  However,
it is still unknown whether $\mathcal{F}^{1}_{(n)} /\mathcal{F}^{1}_{(n.5)}$ has infinite rank for
$n\geq 3$.  

For $m\geq 2$, we will show that for all $n \in \mathbb{N}_{0}$, $\mathcal{F}^{m}_{(n)}
/\mathcal{F}^{m}_{(n+1)}$ contains an infinitely generated subgroup, the subgroup ``generated by boundary links'' defined as follows.
We define the $(n)$-solvable filtration of $\SB(m)$, the subgroup of boundary string links, by 
\begin{equation}\SB\mathcal{F}^{m}_{(n)} \cong \SB(m) \cap \mathcal{F}^{m}_{(n)}.\label{bf}\end{equation}  Then for each $n\geq 0$,
$\SB\mathcal{F}^{m}_{(n)}
/\SB\mathcal{F}^{m}_{(n+1)}$ is a subgroup of $\mathcal{F}^{m}_{(n)}
/\mathcal{F}^{m}_{(n+1)}$.
We will show that the
abelianization of $\SB\mathcal{F}^{m}_{(n)}
/\SB\mathcal{F}^{m}_{(n+1)}$ has infinite rank.  
It would be very interesting to know whether $\SB\mathcal{F}^{m}_{(n)}/\SB\mathcal{F}^{m}_{(n.5)}$ is infinitely generated.  It would be
even more interesting to exhibit non-trivial links in $\SB\mathcal{F}^{m}_{(n.5)}/\SB\mathcal{F}^{m}_{(n+1)}$.
To show that $\SB\mathcal{F}^{m}_{(n)}
/\SB\mathcal{F}^{m}_{(n+1)}$ is infinitely generated, we show that $\rho_{n}$ is
a homomorphism on the subgroup $\SB(m)$ and that $\rho_{n}$ vanishes
on $(n+1)$-solvable links.  We begin with the latter.  

\begin{theorem}\label{rho0}If a $3$-manifold $M$ is $(n)$-solvable then
 for each $(n)$-solution $W$ and $k \leq n$, the inclusion $i:M \rightarrow W$ induces monomorphisms
 \begin{equation}i_{\ast}:  H_{1}( M; \mathcal{K}( \pi_{1}(W)/\pi_{1}(W)^{\sss (k)}_{\sss H}) )\hookrightarrow H_{1}( W; \mathcal{K}( \pi_{1}(W)/\pi_{1}(W)^{\sss (k)}_{\sss H}) )\end{equation}
  and
\begin{equation}i_{\ast}: \frac{\pi_{1}(M)}{\pi_{1}(M)^{(k+1)}_{H}} \hookrightarrow \frac{\pi_{1}(W)}{\pi_{1}(W)^{(k+1)}_{H}};\end{equation}
and  \begin{equation}\rho_{k-1}(M)=0.\end{equation}  Thus, if $L \in \mathcal{F}_{(n)}$ then $\rho_{k-1}(L)=0$ for $k\leq n$.
\end{theorem}
\begin{proof}  Let $W$ be an $(n)$-solution for $M$, $G=\pi_{1}(M)$, and $E=\pi_{1}(W)$.   We will prove the result by induction on $k$.  The result is clearly true for $k=0$ (here $\rho_{-1}(M)=0$ for any $M$).  Assume the result is true for some $k\leq n-1$.

Let $W^{\sss (k+1)}_{\sss H}$ be the regular cover of $W$ corresponding to $E_{\sss H}^{\sss (k+1)}$, the $k^{th}$
term of the torsion-free derived series of $E$.  Since $E^{\sss (n)} \subset E^{\sss (n)}_{\sss H} \subset  E^{\sss (k+1)}_{\sss H}$, $W$ admits 
a ``torsion-free $(k+1)$-Lagrangian'' with ``torsion-free $(k+1)$-duals.''
Specifically, there are intersection and self-intersection forms $\lambda^{\sss H}_{\sss k+1}$ and $\mu_{\sss k+1}^{\sss H}$ on $H_{2}(W_{\sss H}^{\sss (k+1)})$. We project 
the $(n)$-Lagrangian $L$ and $(n)$-duals for $L$ to $H_{2}(W_{\sss H}^{\sss (k+1)})$ to get ``torsion free $(k+1)$-surfaces'' $l_{1}, \dots, l_{g}, d_{1},\dots, d_{g}$ such that
\[\lambda_{k+1}^{H}(l_{i},l_{j})=0 
\text{ and }\lambda^{H}_{k+1}(l_{i},d_{j})=\delta_{ij}\]
for $1 \leq i,j \leq g$.

Recall that for a group $E$, $E_{k+1} = E/E_{\sss H}^{\sss (k+1)}$. 
We will show that $\{l_{1}, \dots, l_{g}, d_{1},\dots, d_{g}\}$ is a $\Z E_{k+1}$-linearly independent set in $H_{2}(W_{\sss H}^{\sss (k+1)})$.  
Suppose 
\[0= \sum_{i=0}^{g} m_{i} l_{i} + n_{i} d_{i}\]
for some $m_{i},n_{i} \in \Z E_{\sss k+1}$.  Then applying $\lambda_{k+1}^{H}(l_{j},-)$ we get
 \[0=\sum_{i=0}^{g} m_{i} \lambda_{k+1}^{H}(l_{j},l_{i}) + n_{i} \lambda_{k+1}^{H}(l_{j},d_{i})=n_{j}\] for $1 \leq  j \leq g$.
 We  now apply $\lambda_{k+1}^{H}(-,d_{j})$ to the new equality $0=\sum_{i=0}^{g} m_{i} l_{i}$ to get $m_{j}=0$ for all $1 \leq j \leq g$.  Therefore $\{l_{i},d_{i} | 1\leq i \leq g\}$
 is a linearly independent set. 
 Thus $\{l_{1}, \dots, l_{g}, d_{1},\dots, d_{g}\}$ generates a rank $2g$ free $\Z E_{\sss k+1}$-submodule of $H_{2}(W_{\sss H}^{\sss (k+1)})$. 
Moreover, since $$\rk_{\mathcal{K}{E_{k+1}}}H_{2}(W;\mathcal{K}E_{k+1})\leq \beta_{2}(W)=2g,$$ by Proposition 4.3 of \cite{COT}, 
$\{l_{1}, \dots, l_{g}, d_{1},\dots, d_{g}\}$ generates $H_{2}(W;\mathcal{K}E_{k+1})$ as a $\mathcal{K}E_{\sss k+1}$-module.   In particular, $H_{2}(W;\mathcal{K}E_{k+1})$ is a free $\mathcal{K}E_{k+1}$-module of rank $2g$.

 Consider the homomorphism $\overline{h}_{W,E_{k+1}}:H_{2}(W;\mathcal{K}E_{k+1})\rightarrow H_{2}(W;\mathcal{K}E_{k+1})^{\ast}$ defined by replacing the coefficients $\mathcal{U}\Lambda$
 in Section~\ref{section:def_rhon} ($\ref{eq:herm}$) with $\mathcal{K}\Lambda$ where $\Lambda=E_{k+1}$.  
 Since $\overline{h}_{W,E_{k+1}}(d_{i})(l_{j})=\overline{\lambda_{k+1}(l_{i},d_{j})}\otimes\mathcal{K}E_{k+1}=\delta_{ij}$ and $\overline{h}_{W,E_{k+1}}(l_{i})(l_{j})=
 \lambda_{k+1}(l_{i},l_{j})\otimes \mathcal{K}E_{k+1}=0$ for $1\leq i,j \leq g$, we can use the same argument as before to show that if 
$\overline{h}_{W,E_{k+1}}(\sum_{i=0}^{g} m_{i} l_{i} + n_{i} d_{i})=0$ then $m_{i}=n_{i}=0$ for all $1\leq i\leq g$. Hence $\overline{h}_{W,E_{k+1}}$ is a monomorphism.  However,
since $\overline{h}_{W,E_{k+1}}$ is a monomorphism between modules of the same rank, $\overline{h}_{W,E_{k+1}}$ is an isomorphism.  In particular, the map 
$H_{2}(W;\mathcal{K}E_{k+1})\rightarrow H_{2}(W,M;\mathcal{K}E_{k+1})$ is surjective which implies that the map $H_{1}(M;\mathcal{K}E_{k+1})\rightarrow H_{1}(W;\mathcal{K}E_{k+1})$
is a monomorphism.  Hence, by Theorem~\ref{CH_thm}, 
\begin{equation}\label{mono}i_{\ast}: \frac{\pi_{1}(M)}{\pi_{1}(M)^{(k+2)}_{H}} \hookrightarrow \frac{\pi_{1}(W)}{\pi_{1}(W)^{(k+2)}_{H}}\end{equation}
 is a monomorphism. Moreover,  
 $(W, \phi_{k} : E \rightarrow E_{k+1} ,i_{\ast}: G_{k+1} \rightarrow E_{k+1})$ is therefore an 
 s-nullbordism for 
 $(M, \phi_{k} :  G \rightarrow G_{k+1} )$.  By definition, $\rho_{k}(M)= \sigma^{(2)}(W,E_{k+1}) -\sigma(W)$.
 
    Since $\sigma(W)=0$, to complete the proof, we show that $\sigma^{(2)}(W,E_{k+1})=0$.
 Recall that $H_{2}(W,M;\mathcal{K}E_{k+1})\cong H_{2}(M;\mathcal{K}E_{k+1})^{*}\cong (\mathcal{K}E_{k+1})^{2g}$ where the isomorphism is given 
 by the composition of the Poincare duality and Kronecker map.
Thus $H_{2}(W;\mathcal{K}E_{k+1})\rightarrow H_{2}(W,M;\mathcal{K}E_{k+1})$ is a surjective map between finitely 
generated $\mathcal{K}E_{k+1}$-modules of the same rank; hence is an isomorphism.   Since $E_{k+1}$ is an Ore domain, $\mathcal{U}E_{k+1}$ is flat as a right 
$\mathcal{K}E_{k+1}$-module hence $H_{2}(W;\mathcal{U}E_{k+1})\rightarrow H_{2}(W,M;\mathcal{U}E_{k+1})$ is an isomorphism.  
 By naturality, the $(k+1)$-Lagrangian above also becomes a metabolizer for $h_{W,E_{k+1}}$ with $\mathcal{U}E_{k+1}$ coefficients. Therefore, $h_{W,E_{k+1}}$ is a nonsingular pairing with 
metabolizer so by Lemma~\ref{witt}, $\sigma^{(2)}(W,E_{k+1})=0$. 
\end{proof}

We now show that $\rho_{n}$ is additive on boundary links.  We begin with an algebraic lemma that will be useful in the proof.

\begin{lemma}\label{retr}Let $A$ be a finitely related group, $E$ be a finitely generated group and $i:A \into E$ be a monomorphism that induces an 
isomorphism on $H_{1}(-;\Q)$.  If there is a retract $r: E \rightarrow A$
then for each $n\geq0$, both $i$ and $r$ induce an isomorphism $i_{\ast} : A/\an \into E/E_{\sss H}^{\sss (n)}$ and $r_{\ast}:E/E_{\sss H}^{\sss (n)} \rightarrow A/\an$ respectively.
Moreover, for each $n\geq 0$, $$\rk_{\mathcal{K}(E/E_{\sss H}^{\sss (n)})} H_{1}(E;\mathcal{K}(E/E_{\sss H}^{\sss (n)}))=
\rk_{\mathcal{K}(A/\an)} H_{1}(A;\mathcal{K}(A/\an)) .$$
\end{lemma}
 \begin{proof}Since $r$ is a retract and $i$ induces an isomorphism on $H_{1}(-;\Q)$, $r$ induces an isomorphism on $H_{1}(-;\Q)$ and a surjective map $H_{2}(-;\Q)$.  
 Hence, by Theorem~4.1 of \cite{CH} (see also Theorem~\ref{CH_thm} and Proposition~\ref{CHprop} of this paper),  $r$ induces isomorphisms $r_{\ast}:E/E_{\sss H}^{\sss (n)} \rightarrow A/\an$ for all $n\geq 0$.
 We will prove $i$ induces monomorphisms $i_{\ast} : A/\an \into E/E_{\sss H}^{\sss (n)}$ by induction.  This is clear when $n=0,1$.  Assume that the result holds for some $n\geq 1$, then by Proposition~\ref{functoriality}, $i$ induces a homomorphism
 $i_{\ast}^{n+1} : A/\anp \rightarrow E/E_{\sss H}^{\sss (n+1)}$.   Moreover, if we postcompose $i_{\ast}^{n+1}$ with $r_{\ast}:E/E_{\sss H}^{\sss (n+1)} \xrightarrow{\cong} A/\anp$, we get the identity.  Hence $i_{\ast}^{n+1}$ is an isomorphism. 
The last statement follows from the last part of Theorem~\ref{CH_thm}. 
\end{proof}

\begin{figure}[htb]
\begin{picture}(44,96)
\put(0,0){{\includegraphics{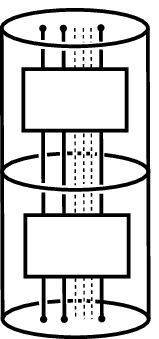}}}
\put(20,67){$L_{1}$}
\put(20,26){$L_{2}$}
\put(-35,49){$L_{1}L_{2}=$}
\end{picture}
\caption{The product of $L_{1}$ and $L_{2}$, $L_{1}L_{2}$}\label{fig:stack}
\end{figure}

If $L_{1}$ and $L_{2}$ are $m$-component string links then $L_{1}L_{2}$ is the $m$-component string link obtained by stacking $L_{1}$ on top of $L_{2}$ as depicted in Figure~\ref{fig:stack}.
This stacking operation induces the multiplication in the group $\mathcal{C}(m)$.

\begin{proposition}\label{rhoadd} If ${L}_1, {L}_2 \in \mathcal{C}(m)$  and $\widehat{L}_1, \widehat{L}_2 $ are boundary links then for each $n\geq 0$,
 $$\rho_n(L_1 L_2)=\rho_n(L_1)+\rho_{n}(L_2).$$
\end{proposition}
\begin{proof}  
If $\widehat{L}_{1}$ and $\widehat{L}_{2}$ are boundary links then $\widehat{L_{1}L_{2}}$ is also a boundary link.
Let $M_{1}$, $M_{2}$ and $M$ be the closed $3$-manifolds obtained by performing $0$-framed Dehn surgery on 
$\widehat{L}_{1}$, $\widehat{L}_{2}$, and $\widehat{L_{1}L_{2}}$ respectively.  Let $M^{\prime}$ be the $3$-manifold obtained by 
performing $0$-framed Dehn surgery along the curves $\alpha_{1},\dots,\alpha_{n}$ in $M_{1}\# M_{2}$ as in Figure~\ref{addhandles}. We will show that 
$M^{\prime}$ is homeomorphic to $M$. 
\begin{figure}[htbp]
\begin{center}
\begin{picture}(180,170)
\put(0,0){\includegraphics{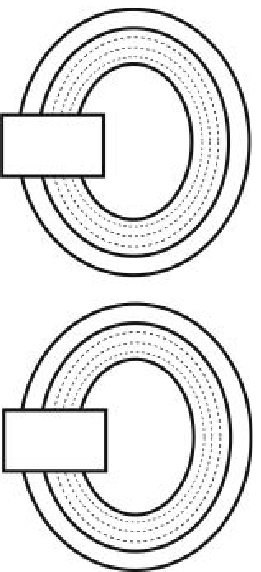}}
\put(10,34){$L_2$}
\put(10,121){$L_1$}
\put(52,37){${\sss 0}$}
\put(69,37){${\sss 0}$}
\put(75,37){${\sss 0}$}
\put(52,124){${\sss 0}$}
\put(69,124){${\sss 0}$}
\put(75,124){${\sss 0}$}
\put(80,81){$\Longrightarrow$}
\put(105,0){\includegraphics{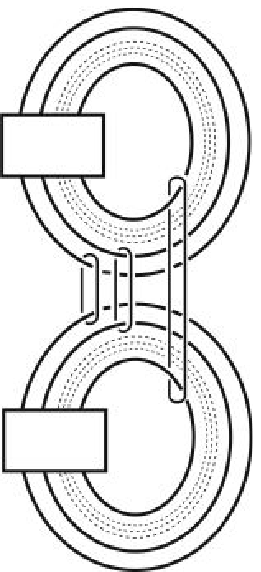}}
\put(115,34){$L_2$}
\put(115,120){$L_1$}
\put(157,37){${\sss 0}$}
\put(174,37){${\sss 0}$}
\put(180,37){${\sss 0}$}
\put(157,124){${\sss 0}$}
\put(174,124){${\sss 0}$}
\put(180,124){${\sss 0}$}
\put(125,81){${\sss 0}$}
\put(135,81){${\sss 0}$}
\put(160,81){${\sss 0}$}
\put(146,81){$\scriptstyle{\cdots}$}
\put(121,77){${\sss {\alpha_{1}}}$}
\put(160,77){${\sss {\alpha_{m}}}$}
\put(112,62){${\scriptsize c}$}\end{picture}
\end{center}
\caption{The $3$-manifolds $M_{1}\# M_{2}$ and $M^{\prime}$}
\label{addhandles}
\end{figure}
To see this, we first isotope the curve $c$ in $M^{\prime}$ (by a ``handle slide'') to obtain the $3$-manifold on the left-hand side of Figure~\ref{slamdunk}. 
Note that by $\widetilde{L}_{i}$, we mean the $(m+1)$-component string link obtained from $L_{i}$ by adding a copy of the ``first string'' by doing a $0$-framed pushoff. 
 Then by doing a slam-dunk \cite[The Slam-Dunk Theorem, p.15]{CG} move on the middle surgery diagram of Figure~\ref{slamdunk}, we see that $M^{\prime}$ is homeomorphic to the $3$-manifold
on the right-hand side of Figure~\ref{slamdunk}.  
\begin{figure}[htbp]
\begin{center}
\begin{picture}(315,165)
\put(0,0){\includegraphics{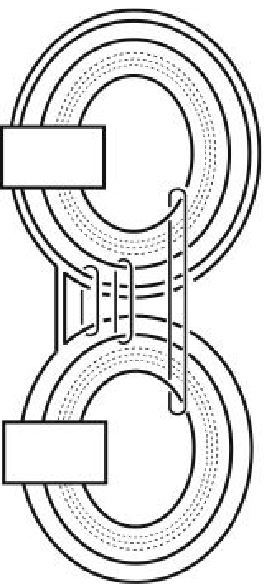}}
\put(10,34){$L_2$}
\put(10,120){$\widetilde{L}_1$}
\put(52,37){${\sss 0}$}
\put(69,37){${\sss 0}$}
\put(75,37){${\sss 0}$}
\put(52,124){${\sss 0}$}
\put(69,119){${\sss 0}$}
\put(67,124){${\sss 0}$}
\put(20,78){${\sss 0}$}
\put(30,78){${\sss 0}$}
\put(40,80){${\sss \dots}$}
\put(55,78){${\sss 0}$}
\put(7,62){${\scriptsize c}$}
\put(80,75){$\Longrightarrow$ }
\put(83,67){$\cong$}
\put(104,0){\includegraphics{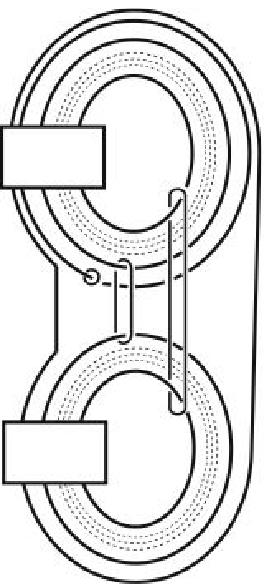}}
\put(114,34){$L_2$}
\put(114,120){$\widetilde{L}_1$}
\put(156,37){${\sss 0}$}
\put(173,37){${\sss 0}$}
\put(179,37){${\sss 0}$}
\put(156,124){${\sss 0}$}
\put(172,124){${\sss 0}$}
\put(172,100){${\sss 0}$}
\put(134,78){${\sss 0}$}
\put(144,80){${\sss \dots}$}
\put(159,78){${\sss 0}$}
\put(125,85){${\sss 0}$}
\put(195,75){$\Longrightarrow$}
\put(199,67){$\cong$}
\put(219,0){\includegraphics{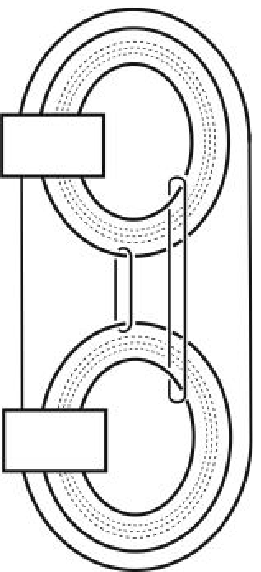}}
\put(230,34){$L_2$}
\put(230,120){$L_1$}
\put(271,37){${\sss 0}$}
\put(288,37){${\sss 0}$}
\put(294,75){${\sss 0}$}
\put(271,124){${\sss 0}$}
\put(287,124){${\sss 0}$}
\put(249,78){${\sss 0}$}
\put(259,80){${\sss \dots}$}
\put(274,78){${\sss 0}$}
\put(183,82){\footnotesize slam-dunk}
\end{picture}
\end{center}
\caption{}\label{slamdunk}
\end{figure}
We continue this process until we arrive at the surgery description of $M^{\prime}$ in Figure~\ref{finalbraid}. Hence $M^{\prime}$ is homeomorphic to $M$.  
\begin{figure}[htbp]
\begin{picture}(76,126)
\put(0,0){\includegraphics{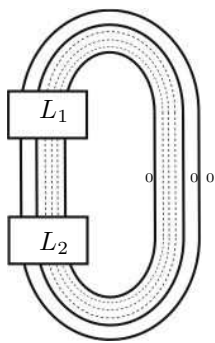}}
\put(12,34){$L_{2}$}
\put(12,84){$L_{1}$}
\put(52,60){${\sss 0}$}
\put(69,60){${\sss 0}$}
\put(75,60){${\sss 0}$}
\end{picture}
\caption{The $3$-manifold $M$, that is homeomorphic to $M^{\prime}$}\label{finalbraid}
\end{figure}

Thus $(M_{1}\# M_{2}) \sqcup \overline{M}=\partial W$ where $W$ is a $4$-dimensional manifold that is obtained by adding $0$-framed
 2-handles to $(M_{1} \# M_{2}) \times I$ along the curves $\alpha_{1}\times \{1\}, \dots, \alpha_{m} \times \{1\}$. 
Let $E=\pi_{1}(W)$, $A=\pi_{1}(M_{1})$, $B=\pi_{1}(M_{2})$ and $\mu^{j}_{1},\dots \mu^{j}_{m}$  be the standard meridians of $L_{j}$ included into $M_{j}$.
Then 
$$E\cong A*B/\left<\mu_{1}^{1}=\mu_{1}^{2},\dots,\mu_{m}^{1}=\mu_{m}^{2}\right>.$$
Recall that 
$M_1 \sqcup M_2 \sqcup \overline{M_1 \# M_2}=\partial W^{\prime}$ where $W^{\prime}$ is obtained by adding a 1-handle
to $(M_1 \sqcup M_2)\times I$.  Let $V=W \cup_{M_{1}\# M_{2}} W^{\prime}$ then $\partial V=M_{1} \sqcup M_{2} \sqcup \overline{M}$ and the inclusion map of 
W into V induces an isomorphism $\pi_{1}(V) \cong \pi_{1}(W)$.
Let $i,j$ and $k$ be the inclusion maps of $M_{1},M_{2}$ and $M$ into $V$ respectively and let $\Gamma_{n}=E/E_{\sss H}^{\sss (n+1)}$.
To complete the proof we will show that (1) $i_{\ast}:A\rightarrow E$ (similarly for $j,k$) induces an isomorphism $A/A_{\sss H}^{\sss (n+1)} \rightarrow \Gamma_{n}$  and 
(2) $\sigma(V,\Gamma_{n})-\sigma(V)=0$ for each $n\geq 0$.  Thus, $\rho_{n}(M_{1})+\rho_{n}(M_{2})=\rho_{n}(M)$. % by Lemmas~\ref{} and \ref{}.
 
 We first prove (1).  Let $F=<x_{1}, \dots, x_{m}>$ be the free group on $m$ generators and $i_{B}: F \rightarrow B$ and $i_{A}: F \rightarrow A$ be the 
 inclusion maps of sending $x_{i} \mapsto \mu_{i}^{j}$. Then $E \cong A *_{F} B$.  Since $B$ is the fundamental group of 0-surgery on a boundary link, 
 there is a retract $r_{B} : B \rightarrow F$ giving a retract $r: E \rightarrow A$ where the inclusion is $i_{\ast}$ as above. Therefore, by Lemma~\ref{retr}, 
 $A/\anp \rightarrow \Gamma_{n}$ is a monomorphism for all $n\geq 0$.  The proof for $B$ is the same as for $A$. 
 
Let $G=\pi_{1}(M)$.  We will show that $k_{\ast}$ induces an isomorphism $G/\gn \rightarrow E/E^{\sss (n)}_{\sss H}$ for each $n\geq 0$.
Since $A$ has a retract to $F$ and $i_{A}$ induces an isomorphism on $H_{1}$, by Lemma~\ref{retr}, $\rk_{\mathcal{K}A_{n}}H_{1}(A;\mathcal{K}A_{n})=
\rk_{\mathcal{K}F_{n}}H_{1}(F;\mathcal{K}F_n)=m-1$ for all $n\geq 0$.  Similarly, $\rk_{\mathcal{K}A_{n}}H_{1}(A;\mathcal{K}A_{n})=
\rk_{\mathcal{K}E_{n}}H_{1}(E;\mathcal{K}E_n)$.  Hence, for all $n\geq 0$, $\rk_{\mathcal{K}\Gamma_{n}}H_{1}(E;\mathcal{K}\Gamma_{n})=
m-1$.  Since $G$ is also the fundamental group of 0-surgery on a boundary link, $\rk_{\mathcal{K}G_{n}}H_{1}(G;\mathcal{K}G_{n})=m-1$. 

We will show by induction on $n$ that $k_{\ast}:G\rightarrow E$ induces an isomorphism $H_{1}(G;\mathcal{K}E_{n}) \rightarrow H_{1}(E;\mathcal{K}E_{n})$ for all $n\geq0$.  
First, $E=G/<[\mu_{i},g_{i}]>$ where $\mu_{i}$ are the given meridians for the string link
$L_{1}L_{2}$ included into $G$.  Hence $H_{1}(G;\mathcal{K}E_{0}) \rightarrow H_{1}(E;\mathcal{K}E_{0})$ is an isomorphism.   
Now, assume that $H_{1}(G;\mathcal{K}E_{n}) \rightarrow H_{1}(E;\mathcal{K}E_{n})$ is an isomorphism for all $k\leq n$.  By Theorem~\ref{CH_thm},  
$G/\gnp \xrightarrow{\cong} E/E^{\sss (n+1)}_{\sss H}$ is an isomorphism.   Hence, $\rk_{\mathcal{K}E_{n+1}}H_{1}(G;\mathcal{K}E_{n+1})=
\rk_{\mathcal{K}G_{n+1}}H_{1}(G;\mathcal{K}G_{n+1})=m-1$.  
Since $k_{\ast}:G \rightarrow E$ is surjective, $k$ induces a surjective map $H_{1}(G;\mathcal{K}E_{n+1}) \onto H_{1}(E;\mathcal{K}E_{n+1})$ between $\mathcal{K}E_{n+1}$-modules 
of the same rank; hence is an isomorphism.  Thus, for each $n\geq0$, $G/\gn \rightarrow E/E_{\sss H}^{\sss (n)}$ is an isomorphism by Theorem~\ref{CH_thm}.

To prove (2), note that there is an exact sequence coming from the long exact sequence of the pair $(W,M_{-})$ with coefficients in $\mathcal{K}(\Gamma_{n})$ 
\[0 \rightarrow \text{Im}(\alpha)  \rightarrow H_{2}(W,M_{-}) \rightarrow H_{1}(M_{-}) \rightarrow H_{1}(W) \rightarrow H_{1}(W,M_{-})\rightarrow 0\]
where $M_{-}=M_{A}\cup M_{B}$ and $\alpha: H_{2}(W) \rightarrow H_{2}(W,M_{A}\cup M_{B}) $.  Since $W$ is obtained by adding a $1$-handle and $m$ $2$-handles to $M_{-}\times I$,
$\rk_{\mathcal{K}\Gamma}H_{2}(W,M_{-})-\rk_{\mathcal{K}\Gamma}H_{1}(W,M_{-})=m-1$.  Since $\rk_{\mathcal{K}\Gamma}H_{1}(W)=m-1$ and $\rk_{\mathcal{K}\Gamma}H_{1}(M_{-})=2m-2$,
it follows that $\rk_{\mathcal{K}\Gamma}\text{Im}(\alpha)=0$. Hence $H_{2}(M_{-})\rightarrow H_{2}(W)$ is a surjective homomorphism.  
\end{proof}

\begin{corollary}\label{cor:hom}For each $n\geq 0$ and $m\geq 1$, $\rho_{n}: \mathcal{B}(m) \rightarrow \R$ is a homomorphism.
\end{corollary}

We now establish the main theorems of this paper.

\begin{theorem}\label{main1} For each $n\geq 0$ and $m \geq 2$, the abelianization of $\SB\mathcal{F}_{(n)}^m/\SB\mathcal{F}_{(n+1)}^m$
has infinite rank.  In particular, for $m\geq 2$, $\SB\mathcal{F}_{(n)}^m/\SB\mathcal{F}_{(n+1)}^m$ is an infinitely generated subgroup of $\mathcal{F}_{(n)}^m/\mathcal{F}_{(n+1)}^m$.
\end{theorem}
\begin{proof}
Let $T$ be the m-component trivial  link in $S^{3}$ with $m\geq 2$ then $F = \pi(S^{3}-T)=\pi_{1}(M_{T})$ is free group on $m$-genererators.  Let $\eta$ be a curve in $S^{3}-T$ such that $\eta$ is a trivial knot and 
the homotopy class of $\eta$ is in $F^{(n)}-F^{(n+1)}$.  Define 
\begin{equation}\label{set_eta}  S_{\eta} = \{T(\eta,K) \hspace{2pt}| \hspace{2pt} K \text{ is a knot in }S^{3}\text{ with Arf invariant zero} \} \subset \SB(m).\end{equation}  
Since the trivial link is $(n)$-solvable,
by the proof Proposition~3.1 of \cite{COT1} (this result holds if one replaces $(n)$-solvable knot with $(n)$-solvable link), $T(\eta,K)$ is $(n)$-solvable if the Arf invariant of $K$ is $0$. Hence $S_{\eta} \subset \SB\mathcal{F}_{(n)}^{m}$.
By Theorem~\ref{rho0}, $\rho_{n}$ vanishes on $\mathcal{F}_{(n+1)}^{m}$. Therefore 
$$\rho_{n}: \frac{\SB\mathcal{F}_{(n)}^{m}}{\SB\mathcal{F}_{(n+1)}^{m}} \rightarrow \R$$
is a homomorphism.

By Proposition~\ref{diff_rho}, $\rho_{n}(T(\eta,K))=\rho_{0}(K)$.  Moreover, by Proposition~2.6 of \cite{COT1}, there is an infinite set of Arf invariant zero knots $\{J_{i}\}$ such that $\{\rho_{0}(J_{i})\}$ is 
$\Z$-linearly independent.  Since $\R$ is abelian, this implies that the abelianization of ${\SB\mathcal{F}_{(n)}^{m}}/{\SB\mathcal{F}_{(n+1)}^{m}} $ has infinite rank.
\end{proof}

Since adding a local knot to $L\in \SB(m)$ doesn't change $\rho_{n}-\rho_{0}$, we can show that $\SB\mathcal{F}_{(n)}^m/\SB\mathcal{F}_{(n+1)}^m$ is infinitely generated ``modulo local knotting.''  
We make this precise starting with the following definition.  Let $K(m)$ be the subgroup of $\SB(m)$ of split string links. $K(m)$ is a normal subgroup of $C(m)$ and hence is a normal subgroup of $\SB(m)$.  For each $n\geq 0$,
define $K\mathcal{F}_{(n)}^{m}=\mathcal{F}^{m}_{(n)} \cap K(m)$.  Then $K\mathcal{F}_{(n+1)}^{m}$ is a normal subgroup of $K\mathcal{F}_{(n)}^{m}$ and 
${K\mathcal{F}_{(n)}^{m}}/{K\mathcal{F}_{(n)}^{m}}$ is a normal subgroup of $\SB\mathcal{F}_{(n)}^m/\SB\mathcal{F}_{(n+1)}^m$.  Note that adding a local knot to $L\in \SB(m)$
is the same as multiplying $L$ by an element of $K(m)$.
As a corollary of the proof of Theorem~\ref{main1} we have the following result.

\begin{corollary}[Theorem~\ref{main1} remains true modulo local knotting] \label{cor:localknotting}For each $n\geq 1$ and $m \geq 2$, the abelianization of 
$$\frac{\SB\mathcal{F}_{(n)}^m/\SB\mathcal{F}_{(n+1)}^m}{K\mathcal{F}_{(n)}^{m}/K\mathcal{F}_{(n)}^{m}}$$
has infinite rank; hence $(\SB\mathcal{F}_{(n)}^m/\SB\mathcal{F}_{(n+1)}^m)/(K\mathcal{F}_{(n)}^{m}/{K\mathcal{F}_{(n)}^{m}})$ is an infinitely generated subgroup of 
$(\mathcal{F}_{(n)}^m/\mathcal{F}_{(n+1)}^m)/(K\mathcal{F}_{(n)}^{m}/{K\mathcal{F}_{(n)}^{m}})$.
\end{corollary}
\begin{proof}This follows from the proof of Theorem~\ref{main1} once we show that $\rho_{n}$ vanishes for $(n)$-solvable string links with $n\geq 1$.  To see this, let $L\in K\mathcal{F}^{m}_{(n)}$ where
$n\geq 1$.  Since $L$ is $(n)$-solvable for $n\geq 1$, $\rho_{0}(L)=0$.  Since $L\in K(m)$, $L$ can be obtained from the trivial string link by tying local knots into the strings.  That is, 
$\hat{L}=L_{m}(\eta_{m},K_{m})$ where $L_{i}$ is defined inductively by: $L_{0}$ is the trivial link and $L_{i+1}=L_{i}(\eta_{i},K_{i})$ for some 
 $\eta_{i}\in \pi_{1}(S^{3}-L_{i})-\pi_{1}(S^{3}-L_{i})_{r}^{\sss (1)}$
and knot $K_{i}$.  Hence by applying Proposition~\ref{diff_rho} multiple times we have $\rho_{n}(L)=\sum_{i=1}^{m}\rho_{0}(K_{i})=\rho_{0}(L)=0$.
\end{proof}

\begin{question} Is $\SB\mathcal{F}^{m}_{(n)}/\SB\mathcal{F}^{m}_{(n.5)}$ (modulo local knotting for $m\geq 2$ and $n\geq 1$) infinitely generated for each $n\geq 0$ and $m\geq 1$?
Is $\SB\mathcal{F}^{m}_{(n.5)}/\SB\mathcal{F}^{m}_{(n+1)}$ (modulo local knotting for $m\geq 2$ and $n\geq 1$) infinitely generated for each $n\geq 0$ and $m\geq 1$?
\end{question}

\subsection{Grope Filtration of $C(m)$}
There is another filtration of the link concordance group called the \emph{grope filtration}, $\mathcal{G}^{m}_{n}$, that is more geometric than the $(n)$-solvable filtration.  

\begin{definition} A grope is a special pair (2-complex, base circles). A grope 
has a height $n \in \mathbb{N}$. A \textbf{grope of height $1$} is precisely a compact, oriented surface $\Sigma$ 
with a non-empty boundary. A \textbf{grope of height  $(n + 1)$} is defined inductively as follows: Let $\{\alpha_{i} , i = 1, \dots ,2g\}$ be a standard 
symplectic basis of circles for $\Sigma$, the bottom stage of the grope. Then a grope 
of height $(n + 1)$ is formed by attaching gropes of height n (with a single boundary component, called the base circle) to each $\alpha_{i}$ along 
the base circle.  A model of a grope can be constructed in $\R^{3}$ and thus has  a regular neighborhood.  Viewing $\R^{3}$ as $\R^{3}\times \{0\} \hookrightarrow \R^{3}\times [-1,1]$, this model grope has a $4$-dimensional regular neighborhood.  
When we say that a grope is embedded in a $4$-dimensional manifold, we always mean that there is an embedding of the entire $4$-dimensional regular neighborhood. 
\end{definition}

We say that $L\in C(m)$ is in $\mathcal{G}^{m}_{n}$ if $\hat{L} \in S^{3}=\partial(D^{4})$ bounds a embedded grope of height $n$ in $D^{4}$.  
Note that if $L_{1},L_{2}\in C(m)$ and $\widehat{L_{i}}$ for $i=1,2$ bounds an embedded grope of height $n$ in $D^{4}$ then $\widehat{L_{1}L_{2}}$ bounds an embedded grope of  height $n$ in $D^{4}$.  
Therefore, $\mathcal{G}^{m}_{n}$ is a subgroup 
of $C(m)$.  Moreover, since $\widehat{L_{1}L_{2}L_{1}^{-1}}=\widehat{L_{2}}$, $\mathcal{G}^{m}_{n}$ is a normal subgroup of $C(m)$. We call 
$$0 \subset \cdots  \subset \mathcal{G}^{m}_{n} \subset \cdots \subset \mathcal{G}^{m}_{1} \subset C(m)$$
the \textbf{grope filtration} of $C(m)$.  
We define the grope filtration of the concordance group of boundary string links by $\SB\mathcal{G}_{n}^{m} = \mathcal{G}_{n}^{m} \cap \SB(m)$.  For more about the grope filtration of a knot, see \cite{COT,COT1,CT}.   For more about gropes, see \cite{FQ,FT}.

The $(n)$-solvable and grope filtrations are related by the following theorem of T. Cochran, K. Orr, and P. Teichner.

\begin{theorem}[Theorem 8.11 of \cite{COT}] If a link $L$ bounds a grope of height $(n+2)$ in $D^{4}$ then $L$ is $(n)$-solvable.
\end{theorem}

Hence for all $n\geq 0$ and $m\geq 1$, $\mathcal{G}_{n+2}^{m} \subset \mathcal{F}_{(n)}^{m}$ (and hence $\SB\mathcal{G}_{n+2}^{m} \subset \SB\mathcal{F}_{(n)}^{m}$).  Note that Cochran-Orr-Teichner only state the above theorem for knots but their proof holds for links in $S^{3}$ as well. 
We show that certain quotients of the grope filtration are non-trivial.  

 \begin{theorem}\label{main2}For each $n\geq 1$ and $m \geq 2$, the abelianization of 
$\SB\mathcal{G}_{n}^m/\SB\mathcal{G}_{n+2}^m$
has infinite rank; hence $\SB\mathcal{G}_{n}^m/\SB\mathcal{G}_{n+2}^m$ is an infinitely generated subgroup of $\mathcal{G}_{n}^m/\mathcal{G}_{n+2}^m$.
\end{theorem} 
\begin{proof} Let $\eta$ be a curve such that the homotopy class of $\eta$ is in $F^{(n)}-F^{(n+1)}$ where $F=\pi_{1}(S^{3}-T)$ and $n \geq 0$.   
By Lemma~3.9 of \cite{CT}, after changing $\eta$ by a homotopy in $S^{3}-T$, we 
can assume that $\eta$ bounds a disk in $S^{3}$ and that $\eta$ bounds an embedded height $n$ grope in $S^{3}-\text{N}(T)$.
Consider the set $S_{\eta}$ as defined in (\ref{set_eta}) (but now using these specially chosen isotopy classes for $\eta$).  We showed in the proof of Theorem~\ref{main1} that $\rho_{n}(S_{\eta})$ is a $\Z$-linearly independent subset of $\R$.
We will show that if $L\in S_{\eta}$ then $L$ bounds a grope of height $(n+1)$.   Since $\SB\mathcal{G}_{n+3}^{m} \subset \SB\mathcal{F}_{(n+1)}^{m}$, 
$\rho_{n}: \SB\mathcal{G}_{n+1}^{m}/\SB\mathcal{G}_{n+3}^{m} \rightarrow \R$ is a well defined homomorphism.  This will complete the proof of the theorem.  

Let $L \in S_{\eta}$ then $L=T(\eta,K)$ where $K$ is a knot.  By Murakami and Nakanishi \cite{MN}, $K$ can be obtained from the unknot by doing a sequence of delta moves.  Hence, in the language of Habiro \cite{Hab}, 
$K$ is related to the unknot by a finite sequence of simple $C_{2}$-moves and ambient isotopies.  Therefore, by Theorem 3.17 of \cite{Hab}, $K$ is the result of clasper surgery on the unknot along $\sqcup_{i=1}^{l} C_{(T_{i},r_{i})}$
where $(T_{i},r_{i})$ is a rooted symmetric tree of height $1$ and the leaves are copies of the meridian of the unknot.  Therefore, $L=T(\eta,K)$ is the result of clasper surgery on the trivial link $T$ along
 $\sqcup_{i=1}^{l} C_{(T_{i},r_{i})}$ where $(T_{i},r_{i})$ is a rooted symmetric tree of height $1$ and the leaves are copies of $\eta$. Since $\eta$ bounds an embedded height $n$ grope in $S^{3}-\text{N}(T)$, 
 by Corollary~3.14 of \cite{CT},  $L$ bounds a height $(n+1)$ grope in $D^{4}$. 
\end{proof}

We remark that there are knots that are the result of a union of clasper surgeries on the unknot along rooted trees of height $2$ and have non-zero $\rho_{0}$ (see for example Figure~3.6 of \cite{CT}).  
By the same argument that was used in the proof of Theorem~\ref{main2}, 
we can show that if you choose such a $K$ then  $L=T(\eta,K)$ bounds a grope of height $n+2$ for $n \geq 0$.  As a result, we see that the groups $\SB\mathcal{G}_{n+2}^{m}/\SB\mathcal{G}_{n+3}^{m}$ are non-trivial for $n\geq 0$.

\begin{proposition} \label{prop:614} For each $n\geq 2$ and $m \geq 2$,  the rank of the abelianization of 
$\SB\mathcal{G}_{n}^m/\SB\mathcal{G}_{n+1}^m$ is at least $1$; hence $\mathcal{G}_{n}^m/\mathcal{G}_{n+1}^m$ contains an infinite cyclic subgroup.
\end{proposition}

Just as in the case of the $(n)$-solvable filtration, Theorem~\ref{main2} and Proposition~\ref{prop:614} are true ``modulo local knotting.''  We formalize this below. Thus these results cannot be obtained using the results 
in \cite{COT,COT1,CT}.  
For each $n\geq 0$,  define $K\mathcal{G}_{n}^{m}=\mathcal{G}^{m}_{n} \cap K(m)$.  Then $K\mathcal{G}_{n+1}^{m}$ is a normal subgroup of $K\mathcal{G}_{n}^{m}$ and 
${K\mathcal{G}_{n}^{m}}/{K\mathcal{G}_{n}^{m}}$ is a normal subgroup of $\SB\mathcal{G}_{n}^m/\SB\mathcal{G}_{n+1}^m$. The proof of the following corollary is similar to the proof of Corollary~\ref{cor:localknotting} so we will omit the proof.

\begin{corollary} $(1)$ For each $n\geq 3$ and $m \geq 2$, the abelianization of 
$$\frac{\SB\mathcal{G}_{n}^m/\SB\mathcal{G}_{n+2}^m}{K \mathcal{G}_{n}^m / K\mathcal{G}_{n+2}^m}$$
has infinite rank;  hence ${(\SB\mathcal{G}_{n}^m/\SB\mathcal{G}_{n+2}^m)}/{(K\mathcal{G}_{n}^m / K\mathcal{G}_{n+2}^m)}$ is an infinitely generated subgroup of 
${(\mathcal{G}_{n}^m/\mathcal{G}_{n+2}^m)}{(K\mathcal{G}_{n}^m/K\mathcal{G}_{n+2}^m)}$.\\ 

\noindent  $(2)$ For each $n\geq 3$ and $m \geq 2$,  the rank of the abelianization of 
$${(\SB\mathcal{G}_{n}^m/\SB\mathcal{G}_{n+1}^m)}{(K\mathcal{G}_{n}^m / K\mathcal{G}_{n+2}^m)}$$ is at least $1$; 
hence $(\mathcal{G}_{n}^m/\mathcal{G}_{n+1}^m)/(K\mathcal{G}_{n}^m / K\mathcal{G}_{n+2}^m)$ contains an infinite cyclic subgroup.
\end{corollary}

\section{Applications to Boundary Link Concordance}\label{sec:last}

We point out some applications of our work to the study of \emph{boundary link concordance} and to the abstract determination of certain $\G$-groups and relative $L$-groups that have been previously 
studied by Cappell-Shaneson and  LeDimet in the context of the classification of links up to concordance (sometimes called cobordism).

Recall that boundary links are amenable to classification because each component bounds a disjoint Seifert surface (alternatively because the fundamental groups of their exteriors admit epimorphisms to the free group).
 In fact it was originally hoped that every odd-dimensional link was concordant to a boundary link, so the classification of link concordance would have been reduced to the case of boundary links. Despite the collapse of
  this hope \cite{CO}, boundary links remain an important case for study.

S. Cappell and J. Shaneson first considered pairs $(L,\th)$ called $\textbf{F-links}$ where $L$ is an $m$-component boundary link and $\th$ is a fixed map $\pi_1(S^3\backslash L)\to F$ that is a splitting map for a 
meridional map. They defined a suitable concordance relation, called $\textbf{F-concordance}$, which entailed an ordinary link concordance between $L_0$ and $L_1$ but required that the fundamental group of the exterior 
of the concordance admit a map to $F$ extending $\th_0$ and $\th_1$. See \cite{CS} for details. 
Let the $m$-component $F$-concordance classes of ($n$-dimensional) $F$-links in $S^{n+2}$ be denoted by $\mathbf{CF(n,m)}$. This is an abelian group if $n>1$, or if $n=m=1$ in which case it is equal to the classical knot concordance group. The \textbf{boundary concordance group of boundary links},
$\textbf{B(n,m)}$, is obtained by dividing out be 
the action of $\aut_{0}(F)$, the group of generator-conjugating automorphisms of the free group, which eliminates the dependence on choice of $\th$. 
This classification (for $n>1$) was later accomplished by K. Ko \cite{Ko}
and W. Mio \cite{Mio} in terms of Seifert matrices and Seifert forms. More recently, D. Sheiham completed a more explicit classification in terms of signatures associated to quivers \cite{Sh}.

For a link $L$ in $S^{n+2}$ (with $n>1$), define $M_{L}$ to be the $(n+2)$-dimensional manifold obtained by doing surgery on $S^{n+2}$ along the components of $L$ 
with the unique normal framing.  For each link $L$ in $S^{n+2}$ where $n \equiv 1 \mod 4$ with $n>1$, and each $k\geq 0$, we define $\rho_{k}(L)=\rho_{k}(M_{L})$.
It is then relatively straightforward (see Proposition~\ref{prop:higher}) to show that the $\rho$-invariants considered herein give a rich source of invariants of $CF(n,m)$ ($n \equiv 1 \text{ mod 4}$ and $n>1$). 
One should compare \cite{L} where $\rho$-invariants associated
 to representations into finite unitary groups are used in an analogous fashion.

\begin{proposition} \label{prop:higher}For any $k\geq 0$, $n\equiv 1 \mod 4$ with $n>1$, the invariant $\rho_k$ induces a homomorphism $\tilde{\rho_k}:CF(n,m)\lra\bbr$ that factors through $B(n,m)$.
\end{proposition}\begin{proof} The situation can be summarized in the following diagram where $\SB(n,m)$ is the group of concordance classes of $m$ component, $n$-dimensional boundary disk links in $D^{n+2}$
(sometimes called boundary string links if $n=1$); where 
$\psi$ is the natural lift defined by Levine \cite[Proposition 2.1]{L} (only for $n>1$); and $I$ is the natural forgetful map.
\[
\begin{diagram}
& & & & \mathcal{B}(n,m) &&\\
&&&\ruTo(4,2)^{\psi}&\dTo &&& \\
CF(n,m)&\rOnto&B(n,m)=CF(n,m)/\aut_0(F)&\rOnto^{I}&\frac{\text{\{Boundary Links\}}}{\text{concordance}}& \rTo^{\rho_{k}}& \R \\
\end{diagram}
\]

Since $\rho_{k}$ is an invariant of concordance of links, the horizontal composition, denoted $\tilde{\rho_{k}}$, exists. \end{proof}

The group $CF(n,m)$, $n>1$, has been classified by the aforementioned authors in terms of $\G$-groups, Seifert matrices and quiver signatures. The biggest question in this field is whether or not 
 $I$ is injective. Our results offer further evidence that it is injective by showing that many powerful signature invariants (the $\rho_k$) of boundary links, that a priori are only invariants of $F$-concordance, are actually
  ordinary concordance invariants.

\begin{question} How many of Sheiham's quiver-signatures are captured by information from $\rho_k$?
\end{question}

\noindent If the $\{\rho_k\}$ were strong enough to detect all of Sheiham's signatures then it would follow that the kernel of $I$ is torsion.

Since Cappell and Shaneson have essentially identified $CF(n,m)$, $n>1$, with the quotient of a certain gamma group $\G_{n+3}(\bbz F\to\bbz)$ (relative $L$-group) modulo the image of $L_{n+3}(\bbz F)$
 \cite[Theorem 2, Theorem 4.1]{CS}
we have the following result.

\begin{proposition} \label{csprop}For each $n \equiv 1 \mod 4$ with $n>1$, and each $k\geq 0$, there is an induced homomorphism
$$
\tilde{\rho}_k: \widetilde\G_{n+3}(\bbz F\to\bbz)/\aut F\lra\bbr.
$$
%that factors through $\SB(n,m)$.
\end{proposition}

Using the techniques of this paper we can also show that each $\tilde{\rho}_k$ extends to the corresponding $\G$-groups of the algebraic closure $\widehat F$ of the free group, in terms of which LeDimet has successfully 
``classified'' the higher-dimensional concordance group of disk links \cite{LD}.

\bibliographystyle{plain}

\bibliography{mybib}

\end{document}